\documentclass[10pt]{amsart}
\usepackage[french, english]{babel}
\usepackage[all]{xy}
\usepackage{graphicx}
\usepackage{psfrag}
\usepackage{amsthm}
\usepackage{amscd}
\usepackage{amsfonts}
\usepackage{amstext}
\usepackage{amssymb}
\usepackage{amsmath}
\usepackage{amsxtra}
\usepackage{plain}

\newcommand{\comm}[1]{}

\numberwithin{equation}{section}
\let\newpf\proof \let\proof\relax 
\newenvironment{pf}{\newpf[\proofname]}{\qed\endtrivlist}

\def\be{\begin{equation}}
\def\ee{\end{equation}}

\def\ba{{\begin{align}}}
\def\ea{{\end{align}}}

\newtheorem{prop}{Proposition}
\newtheorem{cor}[prop]{Corollary}
\newtheorem{ex}{Example}
\newtheorem{lem}[prop]{Lemma}
\newtheorem{thm}[prop]{Theorem}
\newtheorem{defn}{\bf Definition}
\newtheorem{hyp} {Hypothesis}
\newtheorem{rem}{\bf Remark}

\newenvironment{rems}{\noindent {\bf Remarks :} }{}

\newcommand{\RR}{\ensuremath{\mathcal{R}_{b,G}}}
\newcommand{\cl}{\ensuremath{\operatorname{cl}}}
\newcommand{\Ch}{\ensuremath{\operatorname{Ch}}}

\newcommand{\bCh}{\operatorname{\text{}{^b\!}Ch}}
\newcommand{\lb}{\ensuremath{\operatorname{lb}}}
\newcommand{\rr}{\ensuremath{\operatorname{rb}}}
\newcommand{\ff}{\ensuremath{\operatorname{bf}}}
\newcommand{\CCl}{\operatorname{\mathbb C l}}
\newcommand{\Cliff}{\operatorname{Cliff}}
\newcommand{\Cl}{\ensuremath{\operatorname{Cl}}}
\newcommand{\s}{\ensuremath{\mathbf{s}} }
\newcommand{\z}{\ensuremath{\mathbf{z}} }
\newcommand{\x}{\ensuremath{\mathbf{x}} }
\newcommand{\y}{\ensuremath{\mathbf{y}} }

\newcommand{\vv}{\ensuremath{\mathbf{v}} }
\newcommand{\nab}{\ensuremath{\nabla^\mathcal{E}} }
\newcommand{\tnab}{\ensuremath{\tilde{\nabla}^\mathcal{E}} }
\newcommand{\ob}{\ensuremath{\Omega^*\mathcal{B}} }
\newcommand{\oE}{\ensuremath{\Omega^*\mathcal{E}} }
\newcommand{\D}{\operatorname{D}}
\newcommand{\Dt}{\tilde{\operatorname{D}}}
\newcommand{\A}{\operatorname{\mathbb A}}
\newcommand{\Id}{\operatorname{Id}}
\newcommand{\Bb}{\ensuremath{\operatorname{\mathbb{\tilde B}}} }

\newcommand{\Bbs}{\ensuremath{\operatorname{\mathbb{B}_s}} }
\newcommand{\As}{\operatorname{\mathbb{A}_s}}

\newcommand{\Diff}{\ensuremath{\operatorname{Diff}} }
\newcommand{\nabe}{\ensuremath{\operatorname{\nabla^E}} }

\newcommand{\At}{\operatorname{\tilde{\mathbb{A}}}}
\newcommand{\Ats}{\tilde{\mathbb{A}}_s}

\newcommand{\str}{\ensuremath{\text{}{^b\!}STr} }

\newcommand{\bO}{\ensuremath{\text{}{^b\!}\Omega^{\frac{1}{2}}} }
\newcommand{\clf}{\ensuremath{\operatorname{Cliff}(\text{}{^b\!}T^*F)} }
\newcommand{\E}{\ensuremath{\mathcal{E}} }
\newcommand{\B}{\ensuremath{\mathcal{B}} }
\newcommand{\Q}{\mathbb{Q} }
\newcommand{\Z}{\mathbb{Z} }
\newcommand{\N}{\mathbb{N} }
\newcommand{\C}{\mathbb{C} }
\newcommand{\R}{\mathbb{R} }

\long\def\symbolfootnote[#1]#2{\begingroup%
\def\thefootnote{\fnsymbol{footnote}}\footnote[#1]{#2}\endgroup}

\begin{document}

\title{A higher index theorem for foliated 
manifolds with boundary}

\author{Mostafa ESFAHANI ZADEH}

\symbolfootnote[0]{\emph{2000 Mathematics Subject Classification}. 58J22, 53C12.} 
\symbolfootnote[0]{\textit{Key words and phrases}. Groupoid, Index, b-calculus, 
Foliation.}
\address{Mostafa Esfahani Zadeh\\
Institute for Advanced Studies in Basic Sciences (IASBS), Zanjan-Iran}
\email{esfahani@iasbs.ac.ir}

\begin{abstract}
Following Gorokhovsky and Lott and using an extension of the 
b-pseudodifferential calculus of Melrose,
we give a formula for the Chern character of the
Dirac index class of a longitudinal Dirac type operators on a foliated 
manifold with boundary. For this purpose we use the Bismut local index 
formula in the context of noncommutative geometry. This paper uses heavily the
methods and technical results developed by E.Leichtnam and P.Piazza.
\end{abstract}

\maketitle

\begin{otherlanguage}{french}

\begin{center}
{\bf Un th\'{e}or\`{e}me d'indice sup\'{e}rieur pour des vari\'{e}t\'{e}s feuillet\'{e}es \`{a} bord}
\end{center}

\begin{abstract}
En suivant A. Gorokhovsky et J. Lott et en utilisant une extension du b-calcul
de Melrose, nous donnons une formule pour le caract\'{e}re de Chern de la classe d'indice 
d'une famille longitudinale d'op\'{e}rateurs de type Dirac sur une vari\'{e}t\'{e} 
feuillet\'{e}e \`{a} bord. Pour ce faire, 
nous utilisons la formule locale de l'indice de Bismut dans le context de la 
g\'{e}om\'{e}trie non-commutative. 
Cet article utilise de mani\`{e}re essentielle les m\'{e}thodes et techniques 
d\'{e}velopp\'{e}es par E. Leichtnam et P. Piazza. 
\end{abstract}

\end{otherlanguage}

\section{Introduction}
Since its formulation, the Atiyah-Singer index formula
\cite{AtSi-I,AtSi-III} has been, and still is, the subject of a
large amount of significant and interesting researches. The subject
of this paper is closely related to the Family Index Theorem which
is stated and proved by Atiyah and Singer \cite{AtSi-IV}. Consider a
fibration of even dimensional spin manifolds $F\to M\rightarrow B$.
Let $D$ be a family of Dirac type operators acting on smooth sections 
of fiberwise spin bundle twisted by a hermitian vector bundle $E\to B$.
This family gives rise to a continuous map from $B$ into the space
of Fredholm operators. This last space is the classifying space for
topological $K$-theory. So $D$ determines a class ${\rm ind}(D)$
in $K^0(B)$ and its Chern character $\Ch ({\rm ind}(D))$ determines 
an element in $H^*_{dr}(B)$. Using a connection $\nabla^{T(M|F)}$
for the vertical tangent bundle $M$ and a Clifford hermitian 
connection $\nabla^E$, one can define characteristic forms 
$\hat A(TF,\nabla^{T(M|F)})$
and $\Ch(E,\nabla^E)$ which represent elements in $H^*_{dr}(M)$. 
The family index
formula states the following equality in $H^*_{dr}(B)$
\begin{equation}\label{int1}
\Ch({\rm ind}(D))=\frac{1}{(2\pi i)^{p/2}}\int_{M|F}\hat
A(TF,\nabla^{T(M|F)})\wedge \Ch(E,\nabla^E). 
\end{equation}
Here $p=\dim F$ and the integration is performed along 
the fibers and produces a closed differential form which represents an element 
of $H_{dr}^*(B)$. If $B$ is a single point this
reduces to the Atiyah-Singer formula which has a local proof based
on the McKean-Singer formula \cite{McSi-c} and the short 
time asymptotic of the trace density of the heat operator. 
So it is natural to ask
whether there is a local proof for the family index theorem. J.M. 
Bismut proved this family local index theorem by generalizing the 
Daniel Quillen superconnection formalism to infinite dimension
\cite{Bi1,Qu-s}. The index formula for a family of spin manifolds
with boundary was established by Bismut and Cheeger in \cite{BiCh-r}
and generalized significantly by R. Melrose and P. Piazza in \cite{MePi} 
using the $b$-pseudodifferential calculus of Melrose. The
index formula in this case, with the Atiyah-Patodi-Singer 
boundary condition, is the following equality in $H^*_{dr}(B)$
\begin{equation}
\label{int2}\Ch({\rm ind}(D)=\frac{1}{(2\pi i)^{p/2}}\int\hat
A(TF,\nabla^{T(M|F)})\wedge \Ch(E,\nabla^E)
-\frac{1}{2}\eta(D_0). 
\end{equation}
Here $D_0$ is the boundary part of the Dirac family $D$, which is assumed 
to be invertible, and the  differential form $\eta(D_0)$ 
ia a spectral invariant of the boundary family $D_0$. 
The $b$-calculus approach is particularly suitable to handle
the Atiyah-Patodi-Singer boundary conditions.
\newline The family index theorem was partially generalized by 
Connes \cite{Cos}, to state a numerical index theorem for a 
longitudinal foliation family of Dirac type operators. The
foliation is assumed to carry a holonomy invariant measure. The
underlying space for this theorem is the foliation groupoid $G$
associated to the foliated manifold. The generalization of the
Connes' theorem to foliated manifolds with boundary is due to 
Ramachandran \cite{Rama}. In this generalization the foliation is 
assumed to be transverse to the boundary.  
In the foliation case the groupoid structure
provides a convolution algebra structure on the space of compactly
supported function on $G$. A suitable completion of this algebra
form a $C^*$-algebra $C_r^*(G)$. Later, the higher family index
theorem in $K_0(C_r^*(G))$ was established by Connes and
Skandalis \cite{CoSk}. The $C^*$-algebra of a fibration family is
Morita equivalence to the algebra of continuous functions on base
manifold. So this last index theorem reduces to Atiyah-Singer index
theorem for family in $K$-Theory of base space. On the other hand the
Connes' index theorem may be considered as a generalization of
\eqref{int1} to foliation family at the level of zero-forms. The
foliation index theorem in higher degree level was recently
formulated and established by Gorokhovsky and Lott
\cite{GoLo}. In fact A. Gorokhovsky and Lott have proved a
more general index theorem for a family of $G$-invariant Dirac
operator where $G$ is a smooth foliation groupoid. Our work can be
easily extended to include this more general case in the presence of
boundary. They assume that the foliation carries a holonomy
invariant closed current. This current is used to average the higher
degree transversal differential forms to get a number-valued index.
The novelty in their work is to introduce a generalized Chern-Weil
construction for the Chern character of the Dirac index class. Then they 
use the Bismut local index theorem prove their index formula.
In this paper we follow Gorokhovsky and Lott \cite{GoLo} to state
and to prove a higher index formula for foliated manifolds with
boundary. Of course we assume some hypothesis on the foliation and
its holonomy groupoid (see hypothesis \ref{hyp2}) and an
invertibility condition for the boundary family (see hypothesis
\ref{hyp}). In the appendix we recall briefly the generalized
Chern-Weil construction proposed by Gorokhovsky and Lott in
\cite{GoLo}. In section \ref{section1} we fix the notation and set
up the geometric structures. We investigate also some differential
geometric properties of  $b$-foliations which are required in
sequel. In section \ref{section2} we deal with the algebraic and
geometric structures of the holonomy groupoid and prove in
proposition \ref{li} that the generalized curvature of the Bismut
superconnection is $G$-invariant. This is crucial to prove the
Bismut local formula for $G$-invariant family of Dirac operators. In
section \ref{section3} we state, in a $G$-invariant manner, those
aspects of the Melrose's $b$-calculus which are necessary in the 
forthcoming sections. This calculus is the analytic framework for
our work. In particular the $G$-invariant b-calculus with bounds is 
the carrier for the good $G$-invariant parametrices. In this section we state
the $G$-invariant Bismut $b$-density theorem and prove the defect 
formula for the $b$-trace in proposition \ref{ptrace}. 
In section \ref{section4} we
define the   Chern character and study the long time
behavior of the Chern character in proposition \ref{largs}. The
proof of this proposition as well as that of transgression formula
for Chern character in proposition \ref{variation} are
based on the defect formula for $b$-trace. The convergence of eta
form results from some estimates on the heat kernel for
$\Cl(1)$-superconnection which are taken from \cite{GoLo}. 

In this paper we assume that the leaves and the holonomy group of the 
foliation  $(M,F)$ are of polynomial growth. 
Let $E$ be a longitudinal Clifford bundle on $M$ and let $\D=r^*\Dt$ 
be the $G$-invariant Dirac operator acting on the
smooth sections of the $G$-invariant vector bundle 
$r^*E$, where $r:G\to M$ denotes the target map. 
Let $s\colon G\to M$ denote the source map and for each $x\in M$ 
put $G_x:=s^{-1}(x)$. It turns out that $G_x$ is a manifold with boundary and the restriction of 
$r^*E$ to the boundary of this manifold has a natural 
decomposition $r^*E_{|G_x^\partial}=E_{0x}\oplus E_{0x}$.  
The restriction of the family of Dirac operators to the boundary of $G_x$'s 
defines a family $D_0=(D_{0x})_x$ of Dirac operators acting on sections of 
$r^*E_{0x}\to G_x^\partial$. We assume  
the family $D_0$ is invertible. With this assumption, 
the analytical index $ind(\D)$ can be defined as an element in the 
K-theory of a certain subalgebra $\bar{\mathcal R}_{b,G}$ of the algebra of 
b-pseudodifferential smoothing operators with vanishing 
indicial family. 
Given a horizontal connection $\nu:=T^hM$, we define a generalized Chern 
character of this index class taking its value in 
$C^\infty(M,\text{}{^b\!}\Omega^1(F)\otimes\Lambda^*\nu^*)$. 
Let $\rho:C^\infty(M,\text{}{^b\!}\Omega^1(F)\otimes\Lambda^*\nu^*)\to\C$ 
be the linear functional defined by the pairing with a holonomy invariant current. 
For example the integration of function with 
respect to a holonomy invariant transversal measure provides such a functional on 
horizontal differential $0$-forms on $M$.
Let $\mathbb B_s$ be the rescaled $\Cl(1)$-superconnection acting 
on $C^\infty(G^\partial,r^*E_{|\partial})$.  
With the above hypothesis, the 
following integral is proved to be convergent and defines the eta invariant 
of the family $D_0$
\begin{equation*}
\rho(\eta_0)=\frac{1}{\sqrt{\pi}}\int_0^\infty 
\rho\left(STr_\alpha(\frac{d\mathbb B _s}{ds}
e^{-(\mathbb B _s^2-l_0)})\right) \,ds.
\end{equation*}
Here $\mathbb B_s^2-l_0$ is the $G$-invariant part of $\mathbb B_s^2$, 
i.e. it is the pull back of an operator on $(\partial M,\partial F)$ by $r$.
Define the operation $\{~.~\}_p$ by 
the following relation where the longitudinal bundle 
$\text{}{^b\!}\Omega^1(F)$ of $b$-density of order $1$ is identified 
with the longitudinal bundle $\Lambda^p(\text{}{^b\!}TF)$  
\[\{~.~\}_p:\Lambda^*(\text{}{^b\!}TM)
\rightarrow\Lambda^*(T^{h*}M)\otimes\Lambda^p(\text{}{^b\!}TF)
\stackrel{\simeq}{\rightarrow}\Lambda^*(T^{h*}M)\otimes\text{}{^b\!}\Omega^1(F).\]
Under above conditions, results of sections \ref{section3} and \ref{section4} 
imply together the following index theorem (c.f. theorem \ref{hyo}).

\begin{thm}[b-index theorem for $G$-invariant Dirac
operators] Let $G$ denote the holonomy groupoid associated to the
foliated manifold $(M,F)$. 
Let $p=\dim F$ is an even number and let $E$ be
a longitudinal Clifford bundle over $(M,F)$ with
associated  Dirac operator $\Dt$.
If $\rho$ is a linear functional as explained in above then the following 
index formula holds
\begin{equation*}
\rho(\Ch(ind(\D))= \frac{1}{(2\pi
i)^{p/2}}\int_M\rho(x)\{\operatorname{\hat{A}}(M,F)(x)
\operatorname{Ch}(E/S)(x)\}_p-\frac{1}{2}\rho(\eta_0).
\end{equation*}
\end{thm}
The integrand is the Atiyah-Singer characteristic differential form coming 
out from asymptotic behavior of the heat kernel of Bismut superconnection 
for foliation $(M,F)$. The precise definition is given by relations 
\eqref{akar} and \eqref{twi}.\\
\textbf{Acknowledgment: }The author would like to thank the anonymous referee 
for the detailed reviews. We appreciate the careful examination of our paper and 
thank for several important comments..\\
\textbf{Dedication: }The author would like to dedicate this paper to Jamal Amoo Rezaii, 
a great mathematics teacher in Bijar-Iran.

\section{Geometric setting and notation}\label{section1}
\subsection{Dirac longitudinal structures on foliated manifolds 
with boundary}\label{geometric setting}
Let $(M,F)$ be a smooth foliated compact manifold with boundary with
even dimensional leaves which intersect transversally the boundary
$\partial M$. So there exists a neighborhood $U$ of $\partial M$ and
a diffeomorphism $\phi:U\rightarrow \partial M\times [0,1)$ such
that
\begin{equation}\label{ab}
\phi_*(F_{|U})\simeq (\partial M\times [0,1),\partial F\times[0,1))
\end{equation}
where $\phi_*(F_{|U})$ denotes the push-forward of the foliation
$F_{|U}$ by $\phi$. We denote the coordinate of the interval $[0,1)$
by $t$. We denote by $p$ the dimension of the leaves
and by $q=n-p$ their co-dimension.\\
In what follows we refer to \cite{Me} (See also \cite{MePi} and
\cite{Gr}) for notation and basic concepts from b-calculus and to
\cite{LaMi} for basic facts in spin geometry. We denote by $\text{}{^b\!}T
F$ the longitudinal b-tangent bundle of $M$, i.e, the vector bundle
on $M$ whose restriction to each leaf is the b-tangent bundle of
that leaf. 
A longitudinal exact b-metric on $M$ is a smooth bilinear
form on $\text{}{^b\!}T F$ which falls down to an exact b-metric when
it is restricted to each leaf. We suppose that it takes the
following form in the neighborhood $U$
\begin{equation}\label{abc}
g=\frac{d^2t}{t^2}+g_0
\end{equation}
where $g_0$ is a longitudinal Riemannian metric on $T\partial M$.
The set of the smooth sections of the longitudinal half b-density
bundle equipped with the $C^\infty$-topology is denoted by
$C^\infty(M,\bO)$, whereas $\dot{C}^\infty(M,\bO)$ consists of those
elements in $C^\infty(M,\bO)$ which vanish to all order at
$\partial M$. The dual spaces of these spaces are denoted,
respectively,
by $\dot{C}^{-\infty}(M,\bO)$ and $C^{-\infty}(M,\bO)$.\\
Associated to the
riemannian structure $g$ on  $\text{}{^b\!}T^* F$ there is the longitudinal
b-Clifford bundle on $M$ denoted by \clf  which is, due to the 
parity of the dimension of the leaves, a $\Z_2$-graded 
longitudinal vector bundle (in accordance with \cite[page 4]{MePi}, our
convention for the Clifford relation is v.v=g(v,v)). 
Now let $E\rightarrow M$ be a smooth
bundle on $M$ which is a $\clf$ bundle, so it is $\Z_2$-graded and
\(E=E^+\oplus E^-\) with respect to this grading. Let $E$ be
equipped with a hermitian b-connection $\nabe$ which is compatible
with the Clifford action of \clf. We denote by
\(\Dt:C^\infty(M,E)\rightarrow C^\infty(M,E)\) the longitudinal
b-Dirac operator corresponding to these data. This operator is
grading reversing with respect to the grading $E=E^+\oplus E^-$ and
takes the following form
\begin{gather*}
\Dt=\begin{pmatrix}0&\Dt^-\\\Dt^+&0\end{pmatrix}\quad;\quad
(\Dt^+)^*=\Dt^- .
\end{gather*}
The longitudinal Riemannian metric $g_0$ defines the longitudinal
Clifford bundle on the boundary. We denote this bundle by
$\Cliff(T^*\partial F)$. The following application extends to a
natural algebra morphism on this Clifford bundle
\begin{gather}
T^*(\partial F)\overset{j}\rightarrow\clf|\partial M \\
j(\xi)=i\frac{dt}{t}\cdot\xi\notag
\end{gather}
Identify $E^+|\partial M$ by $E^0$. Then $\sigma^-=cl(i\frac{dt}{t})$ is an 
isomorphism between $E^-_{|\partial M}$ and $E^0$.\\
The bundle $E=E^+\oplus E^-$, restricted to the boundary $\partial M$, 
is also a longitudinal Clifford bundle
over $(\partial M, \partial F)$. So we can consider the boundary 
grading reversing Dirac operator
 $\operatorname{\tilde D}_\partial=\operatorname{\tilde D}^+_\partial
 \oplus\operatorname{\tilde D}^-_\partial$. 
It is easy to verify the relation \(\operatorname{\tilde
D}_\partial^-\circ (\sigma^-)^{-1}= \sigma^-\circ
\operatorname{\tilde D}^+_\partial\). Put \(\operatorname{\tilde
D}_0 =\operatorname{\tilde D}_\partial^-\circ
(\sigma^-)^{-1}=\sigma^-\circ \operatorname{\tilde D}^+_\partial\). 
With respect to the isomorphism 
$\sigma:=id\oplus\sigma^-\colon E_{|\partial M}\to E_0\oplus E_0$, 
the Dirac operator $\Dt$ takes the following form as an operator on 
$C^\infty(\partial
M\times [0,1),E_0\oplus E_0)$ 
\begin{equation}\label{nar}
 \Dt_{|\partial M\times[0,1)}=
\begin{pmatrix}0& -t\partial_t+\operatorname{\tilde D_0}\\
    t\partial_t+\operatorname{\tilde D_0}&0\end{pmatrix}.
\end{equation}

\subsection{Differential geometry of foliated manifolds}\label{folman}
As an additional structure, we suppose that there is given a
horizontal distribution $T^hM$, that is a sub bundle of
$\text{}{^b\!}TM$ such that $\text{}{^b\!}TM\simeq
T^hM\oplus\text{}{^b\!}TF$. The first and second component of a vector
$\mathbf{x}\in \text{}{^b\!}TM$ with respect to this decomposition is
denoted, respectively, by $\mathbf{x}^h$and $\mathbf{x}^v$. In what
follows, for simplifying the notation, we denote sometimes the
horizontal bundle by $\nu$. We suppose that $T_x^hM\subset
T_x\partial M$ for $x\in\partial M$, so the sections of $\nu$ are b-vector fields
and the differential operators defined by them are b-
differential operators. The previous decomposition of the tangent
bundle $\text{}{^b\!}TM$ gives rise to the following isomorphism
between graded algebras
 \begin{equation}\label{ff}
\Omega^*(\text{}{^b\!}T^*M)\cong\Omega^*(T^{h*}M)\otimes\Omega^*(\text{}{^b\!}T^*F).
\end{equation}
Let $\omega\in\Omega^*(\text{}{^b\!}T^*M)$ be of degree $(k,l)$ with
respect to this decomposition. The covariant derivation in the
horizontal direction gives an another differential form
$d^h(\omega)$ of degree $(k+1,l)$. We can define $d^h(\mu)$,
where $\mu$ is a density, by identifying the b-densities with powers
of longitudinal differential forms of degree $p$. In this way we get
the graded derivation
  \(d^h:(M;\Lambda^*\nu^*\otimes \bO)\rightarrow
(M;\Lambda^{*+1}\nu^*\otimes \bO)\) by putting
\[d^h(\omega\otimes\mu)
:=d^h\omega\otimes\mu+(-1)^{\deg(\omega)}\omega\wedge d^h\mu.\]
To the horizontal distribution $T^hM$ is associated its curvature
$T^h$. It is a longitudinal vector-valued horizontal 2-form  
defined by the following formula
 \begin{equation}\label{cur}
T^h(\mathbf{y},\mathbf{z})(x)
=[\mathbf{y}^h,\mathbf{z}^h]^v(x)\in\text{}{^b\!}T_xF.
\end{equation}
The Lie derivation along $T^h$ defines an operator
$\mathcal{L}_h:\Omega^*(M)\to \Omega^{*+2}(M)$. Using this operator
one can define obviously an another operator, denoted by the same
symbol,
\begin{gather*}
\mathcal{L}_h:C^\infty(M; \Lambda^*\nu^*\otimes\bO)\rightarrow
C^\infty(M; \Lambda^{*+2}\nu^*\otimes\bO)
\end{gather*}
The definition of $\mathcal{L}_h$ can be directly extended to
 sections of $E$ with differential form coefficients.

The covariant derivation in the horizontal direction defines the
operator $\nabe:\,C^\infty(M,E)\rightarrow C^\infty (M, \nu^*\otimes
E)$. The definition of this operator can be extended to $C^\infty(M,
\Lambda^*\nu^*\otimes\bO\otimes E)$ by putting
\[\nabe(\omega\otimes\mu\otimes\xi):=
d^h\omega\otimes\mu\otimes\xi+(-1)^{\deg(\omega)}\omega\otimes
d^h\mu\otimes\xi+(-1)^{\deg(\omega)}\omega\otimes\mu\otimes\nabe\xi.\]
This means that \nabe is a connection.

The following relations hold between  various operators defined
above (for proofs see, e.g, \cite[section III.7.$\alpha$]{Co} or
\cite[(3.13)]{BiLo})
\begin{equation}\label{curid}
(d^h)^2=-\mathcal{L}_h\quad;\quad (\nabe)^2=-\mathcal{L}_h.
\end{equation}
Using the tensor product of the longitudinal Clifford connection on
$E$ and the Levi-Civita connection on differential forms, one can
construct the longitudinal Dirac operator acting on distributional
sections of $E$. From now on we suppose that the degree of a
longitudinal density is zero, so it commutes with the Clifford
action of the longitudinal (co)tangent bundle. We denote this
operator again by $\Dt$. Now we are ready to define the Bismut superconnection 
$\At$ by
\begin{gather}
\At:\,C^\infty(M,\Lambda^*\nu^*\otimes\bO\otimes E)\rightarrow
C^\infty(M,\Lambda^*\nu^*\otimes\bO\otimes E)\notag\\
\At=\Dt+\nabe-\frac{1}{4i}\cl(T^h).\label{ccc}
\end{gather}
\begin{rem}\label{ai}
let $\mu'$ denote a longitudinal p-form. We recall from 
\cite[lemma 10.4]{BeGeVe} that
$d^M(\mu')=\kappa\wedge\mu'+i_{T^h}(\mu')$. Here $k$ denotes the mean
curvature associated to the degenerated
metric $g=0\oplus g$ on $\text{}^bTM$. In
particular  $d^h(\mu')=\kappa\wedge\mu'$. So we get the following
relation where $\mu$ is a longitudinal half b-density
\begin{equation}
d^h(\mu)=\frac{1}{2}\kappa\wedge\mu.
\end{equation}
This relation shows that the mean curvature term in the definition of the Bismut
superconnection in \cite[lemma 10.4]{BeGeVe} exists implicitly in \eqref{ccc}
in the action of \nabe on half densities.
\end{rem}

\section{Holonomy groupoid}\label{section2}
\subsection{Geometric structures of the holonomy groupoid }

 In this section we establish some differential geometric
constructions which will be used in the forthcoming sections. The
basic reference for the following matters are \cite{Cos, MoSc}.
Let $G$ be the holonomy groupoid associated to the foliation $(M,F)$
and denote the source and the range maps respectively by $s$ and $r$. 
The groupoid $G$ is a manifold with corner and a 
typical trivializing chart for 
it is $(U_x,[\gamma],U_y)$. Here $x$ and $y$ are two points in $M$ which belong to 
the same leaf of the foliation $(M,F)$ while   
$U_x$ and $U_y$ are  trivializing charts of $M$ around $x$ and $y$ respectively.  
$[\gamma]$ is the holonomy class of a path $\gamma$ between $x$ and $y$ where  
$[\gamma]=[\gamma']$ if the holonomy of $\gamma$ and $\gamma'$ are equal. 
It is clear from this discussion that $G_x:=s^{-1}(x)$ and $G^y:=r^{-1}(y)$ 
are smooth manifolds with boundary and the boundary of $G_x$ is $G_x^\partial$ 
where
\[G_x^\partial:=\{u\in G|s(u)=x\text{ and }r(u)\in\partial M\}.\]
The definition of $G^y_\partial$ is similar and we have 
$\partial G^y=G_\partial^y$. In this paper we use the standard notation of 
the theory of groupoid used, e.g. in \cite{MoSc}. For example 
$G^\partial$ consists of those elements $u\in G$ with 
$r(u)\in\partial M$. Clearly $G^\partial$ is the union of $G_x^\partial$ 
for $x\in M$. The similar comments apply to $G_\partial$. 
The boundary $\partial M$ can be used to define the sub-groupoid 
$G^\partial_\partial$ of $G$ by 
\[G^\partial_\partial:=\{u\in G|s(u)\in\partial M
\text{ and }r(u)\in\partial M\}.\]
The source and range maps of this groupoid are 
the restriction of the source and 
range maps of the groupoid $G$ and take their values in $\partial M$. The units of 
the groupoid $G^\partial_\partial$ are naturally identified with the 
points of $\partial M$.

Take $u\in G$ with $s(u)=x$ and $r(u)=y$. Let
$V$ and $V'$ be local charts for the same leaf around, respectively,
$x$ and $y$. Let $T,\,T'$ be sufficiently small local
transversal passing through $x$ and $y$. We denote by $Hol_u:T\rightarrow
T'$ and $u_*:T_xT\to T_yT'$, respectively, the local holonomy
diffeomorphism induced by $u$ and its derivative. The sets
of the form  $\mathcal{U}=V\times V'\times T\times_{Hol_u}T'$ define
a basis for the topology of $G$. With this topology, $G$ has the
structure of  a smooth manifold with corners which is foliated by 
manifolds with boundary $\{G_x\}_{x\in M}$. Given a connection $T^hM$ for
the foliation $(M,F)$, there is a unique horizontal distribution
$T^HM$ on $G$ defined by
\(T^H_uG=(d_xs)^{-1}T^h_xM\cap(d_ur)^{-1}T^h_yM\). The following
relation holds
\begin{gather}\label{chi}
T_uG=T_uG_x\oplus T_uG^y\oplus T^H_uG.
 \end{gather}
Here the first summand is, the so-called, $G$-invariant part of $TG$.
It is given by $r^*T_yF$ when $r$ is restricted to $G_x$. The direct
sum $T_uG_x\oplus T_uG^y$ is the vertical while $T^H_uG$ is the
horizontal part of $T_uG$. The application $dr\oplus ds$ gives an
isomorphism between the vertical part of $T_uG$ and $T_yF\oplus
T_xF$. Concerning the horizontal part, we have the following
isomorphisms given again by $ds$ and $dr$
\begin{gather}\label{too}
\xymatrix{T_x^hM  & T_u^HG \ar[l]^{\sim}_{ds} \ar[r]^{\sim}_{dr} &
T_y^hM.}
\end{gather}
Here the map from the left end to the right end is given by the holonomy $u_*$.
Therefore a horizontal vector $\mathbf{y}$ in $T_uG$ is determined
uniquely by two vector $\mathbf{y}_x\in T^h_xM$ and $\mathbf{y}_y\in
T^h_yM$ such that $u_*(\mathbf{y}_x)=\mathbf{y}_y$. On the other
hand a vertical vector $\mathbf{x}$ in $T_uG$ is determined uniquely
by two vector $\mathbf{x}_x\in T_xF$ and  $\mathbf{x}_y\in T_yF$. A
smooth function $f$ on $G$ get the form $f_V\times f_{V'}\times
f_T\times f_{T'}$ in the local chart $\mathcal{U}$. The holonomy
invariance condition reads  $f_{T'}\circ Hol_u=f_T$. The following
formulas hold and are usually more suitable to perform local
computations
\begin{equation}\label{nos}
\mathbf{x}.f=\mathbf{x}_x.f_V+\mathbf{x}_y.f_{V'}\quad;\quad
\mathbf{y}.f=\mathbf{y}_x.f_T=\mathbf{y}_y.f_{T'}\,.
\end{equation}
Now let $\mathbf{y}$ and $\mathbf{z}$ be two sections of the
distribution $T^HG$ and let $T^H(\mathbf{y}, \mathbf{z})$ denote the
vertical component of the vector field $[\mathbf{y},\mathbf{z}]$.
$T^H$ is in fact a tensor and is called the curvature of the
connection $T^HG$. From the above discussion and using the relation
(\ref{cur}) it is clear that
\begin{gather}\label{dfo}
T^H(\mathbf{y}, \mathbf{z})=
r^*T^h(\mathbf{y}_y, \mathbf{z}_y)\oplus s^*T^h(\mathbf{y}_x,
\mathbf{z}_x).
\end{gather}
So the $G$-invariant part of $T^H$ is given by $r^*T^h$.

The manifold with corners $G$ is foliated by manifolds
$\{G_x\}_{x\in M}$, this is the foliation that we are interested in. We
denote this foliation by $(G,G_*)$. In the sequel we speak about the
longitudinal structures on $G$ with respect to this foliation. Since
the target map $r_{|G_x}$ is a covering map on its range, we can lift
all longitudinal structures on the foliation $(M,F)$ to the holonomy
foliation $(G,G_*)$.
 Such structures are $G$-invariant with respect to the left action of the
groupoid $G$ on itself. We denote by $\D$ the
 lifting of the Dirac operator $\Dt$. For $0\leq k\leq q$ put
\begin{gather*}
\Omega^k\B:=C_c^\infty(G,\Lambda^kr^*\nu^*\otimes
s^*\text{}{^b\!}\Omega^1)\quad;\quad\ob=\oplus_k\Omega^k\B\\
\Omega^k\mathcal{E}=C_c^\infty(G,\Lambda^ks^*\nu^*\otimes
r^*(\bO\otimes E))\quad;\quad\oE=\oplus_k\Omega^k\E.
\end{gather*}
Now we introduce some differential operators on
these spaces. As in the definition of $d^h$, the covariant
derivative along the horizontal direction $T^HG$ defines the
operator $d^H$ on \ob  and the operator \nab on \oE. Let
$r^*(\omega)\otimes s^*(\mu)$ be an element in $\Omega^k\B$, and let
\y be a horizontal vector field. With the notation of the relations
(\ref{nos}) the following formula holds
\begin{gather}\label{pulinv}
d^H(\y)(r^*(\omega)\otimes s^*(\mu))
=r^*(d^h\y_y)(\omega)\otimes s^*(\mu)+(-1)^{\deg\omega}
r^*(\omega)\otimes s^*(d^h\y_x)(\mu).
\end{gather}
Take an element in \oE of the form $s^*(\omega)\otimes
r^*(\mu\otimes\xi)$, a similar local formula holds for the action of
\nab
\begin{align}\label{khod}
\nab(\y)(s^*(\omega)\otimes r^*(\mu\otimes\xi))
=& s^*(d^h\y_x)(\omega)\otimes r^*(\mu\otimes\xi)\notag\\
&+(-1)^{\deg\omega}s^*(\omega)\otimes r^*(d^h\y_y)(\mu)\otimes r^*(\xi)\notag\\
&+(-1)^{\deg\omega}s^*(\omega)\otimes r^*(\mu)\otimes
r^*\nabe(\y_x)(\xi).
\end{align}
 As in the previous section, the Lie derivative
along the curvature vector $T^H$ (with differential form coefficients) 
defines two operators on \ob and \oE that
we denote by the same symbol $\mathcal{L}_H$. The first formula of
(\ref{nos}) with (\ref{dfo}) give the following local formula
\[\mathcal{L}_H(\mathbf{y}, \mathbf{z})(r^*(\omega)\otimes
s^*(\mu))=r^*\mathcal{L}_h(\y_y,\z_y)(\omega)\otimes s^*(\mu)+
r^*(\omega)\otimes s^*\mathcal{L}_h(\y_x,\z_x)(\mu)\] There is a
similar local formula for the action of $\mathcal{L}_H$ on \oE.
Finally we define the operator $l:\,\Omega^k\E\rightarrow\Omega^{k+2}\E$ by 
\begin{gather}\label{k}
l(\mathbf{y}, \mathbf{z})(s^*(\omega)\otimes
r^*(\mu\otimes\xi))=s^*(\omega)\otimes
r^*\mathcal{L}_h(\y_x,\z_x)(\mu\otimes\xi)
\end{gather}
It is clear that $\mathcal{L}_H-l=r^*(\mathcal{L}_h)$, so
$\mathcal{L}_H-l$ is a longitudinal $G$-invariant differential
operator.
\begin{lem}\label{f3}
The following formulas hold
\begin{enumerate}
\item \(\mathcal{L}_H-l=r^*(\mathcal{L}_h).\)
 \item
\((d^H)^2=-\frac{1}{2}\mathcal{L}_H\quad;\quad (\nab)^2=-\mathcal{L}_H.\)
\item
\(\nab\circ l=l\circ\nab\quad;\quad\D\circ l=l\circ\D.\)
\end{enumerate}
In particular $(d^H)^2\neq0$ if the
 connection $T^hM$ is not integrable.
\begin{pf}
The first part is already proved in the above discussion. The second
part follows from the first part and relations (\ref{curid}).
Comparing with the relations (\ref{curid}), the extra factor
$\frac{1}{2}$ is coming from the action on half densities. To prove
the third part, notice that $\D=r^*\Dt$ so $\D(s^*(\omega)\otimes
r^*(\mu\otimes\xi))=(-1)^{\deg(\omega)}s^*(\omega)\otimes
r^*\Dt(\mu\otimes\xi)$. The commutativity relation 
\(\D\circ l=l\circ\D\) is clear from this relation and local expression
(\ref{k}).Using the relations (\ref{khod}) and (\ref{k}), the third
part reduces to relations \(d^h\circ
\mathcal{L}_h=\mathcal{L}_h\circ d^h\) and \(\nabe\circ
\mathcal{L}_h=\mathcal{L}_h\circ\nabe\) which follows from
(\ref{curid}).
\end{pf}
\end{lem}
For the future use we put a $\Cliff(TG)$-module structure on \oE. At
first we equip the tangent bundle $TG$ with a degenerate Riemannian
structure as follows. From (\ref{chi}) at each point $u\in G$ we
have $T_uG=T_uG_x\oplus T_uG^y\oplus T^HG$. The first summand has
the $G$-invariant Riemannian structure $(r_{|G_x})^*g$. Let the two other
summand have the null metric and equip the total space by the direct
sum Riemannian structure. Now take a section of $TG$ which is of the
form $r^*(\x)+s^*(\x')+\vv$ with respect to the previous direct sum
decomposition. Let $\omega\otimes\mu\otimes\x$ be an element in \oE
and put
\begin{gather}
s^*(\x').(\omega\otimes\mu\otimes\xi)=0\quad;\quad
\vv^*.(\omega\otimes\mu\otimes\xi )
=\vv^*\wedge\omega\otimes\mu\otimes\xi\label{cli1}\\
\x.(\omega\otimes\mu\otimes\xi)=(-1)^{deg(\omega)}\omega\otimes\mu\otimes\cl(\x)
.\xi\label{cli2}
\end{gather}
It is clear that $\x'.\x'=\vv.\vv=0$, and $\x.\x=g(\x,\x)$ so these
relations define a Clifford-module structure on $\oE$. This definition
for the Clifford structure and the relation (\ref{dfo}) imply the
$G$-invariance relation $\cl(T^H)=r^*\cl(T^h)$.\\
Now we are ready to introduce the Bismut superconnection 
\(\A:\,\oE\rightarrow \oE\) by the following formula
\begin{gather*}
\A=\D+\nab-\frac{1}{4i}\cl(T^H).
\end{gather*}
The following property of the Bismut superconnection is crucial for
our purpose
\begin{prop}\label{li}
 The following relation holds,
\begin{equation}\label{invar}
\A^2-l=r^*(\At^2)
\end{equation}
where $\At$ denotes the Bismut b-superconnection on the foliated
 manifold $(M,F)$.
\begin{pf}
>From the definition of the superconnection we have
\begin{align}
\A^2-l=&\D^2+\frac{1}{16}\cl(T^H)^2+\frac{-1}{4i}
\{\D\circ \cl(T^H)+\cl(T^H)\circ\D\}\label{square1}\\
&+\{\nab\D+\D\nab\}+\frac{-1}{4i}\{\nab \cl(T^H)+\cl(T^H)\nab\}\label{sq1}\\
       &+\{(\nab)^2-l\}.\notag
\end{align}
The following similar relation holds for the Bismut superconnection $\At$
\begin{align}
\At^2=&\Dt^2+\frac{1}{16}\cl(T^h)^2+\frac{-1}{4i}
\{\Dt\circ \cl(T^h)+\cl(T^h)\circ\Dt\}\label{square2}\\
&+\{\nabe\Dt+\Dt\nabe\}+\frac{-1}{4i}
\{\nabe \cl(T^h)+\cl(T^h)\nabe\}\label{sq2}\\
      &+(\nabe)^2.\notag
\end{align}
We have already shown the relation $\cl(T^H)=r^*\cl(T^h)$. On the
other hand, from the very definition one has  $\D=r^*\Dt$. Thus
each term in the right hand side of (\ref{square1}) is the
pull-back, by $r^*$, of the corresponding term in (\ref{square2}).
On the other hand, the relations (\ref{curid}), with the first and
second part of the lemma \ref{f3} imply the relation
\((\nab)^2-l=r^*(\nabe)^2\). Now consider the general element
$s^\omega\otimes r^*(\mu\otimes\xi)$ in \oE with $\deg{\omega}=k$.
Using the local expression for \nab one has
\begin{gather*}
\nab\D(s^*(\omega)\otimes r^*(\mu\otimes\xi))=(-1)^ks^*\nabe(\omega)\otimes
r^*\Dt(\mu\otimes\xi)+s^*\omega\otimes r^*\nabe\Dt(\mu\otimes\xi)\\
\D\nab(s^*(\omega)\otimes
r^*(\mu\otimes\xi))=(-1)^{k+1}
s^*\nabe\omega\otimes r^*\Dt(\mu\otimes\xi)+s^*\omega\otimes 
r^* \Dt\nabe(\mu\otimes\xi).
\end{gather*}
Adding these relations gives the desired  $G$-invariance
 relation
 \[\nab\D+\D\nab=r^*(\nabe\Dt+\Dt\nabe).\] A quit
similar computation shows the following $G$-invariance relation and
finishes the proof of the proposition
\[\nab\cl(T^H)+\cl(T^H)\nab=r^*(\nabe\cl(T^h)+\cl(T^h)\nabe\]
\end{pf}
\end{prop}
Essentially the same computations in the proof of the last relation
 prove the vanishing relation
 \(\nabe\cl(T^h)+\cl(T^h)\nabe=0\). Now for $s>0$ define the rescaled Bismut
 superconnection by \(\As=s\,\D+\nab-\frac{1}{4is}\cl(T)\).
This vanishing formula shows that $\As^2-l=s^2\,\D^2+sF+J$, where $F$
and $J$ are $G$-invariant differential operators with differential
form coefficients of positive degrees. In particular these operators
are nilpotent. For each $x\in M$, put $G_x^\partial:=\partial G_x$ and let 
$G^\partial$ be the union of all $G_x^\partial$, equipped with the induced 
topology from $G$. 
We recall the isomorphisms $E_0=E^+_{|\partial M}$ and 
$E_0\simeq E^-_{|\partial M}$. 
Using these isomorphisms and the fact that connection 
$T^hM$ is tangent to the boundary at boundary points, it turns out that  
$\nab$, $T^H$ and $l$ (c.f. \eqref{k}) define differential operators 
on $C_c^\infty(G^\partial,\Lambda^ks^*\nu^*\otimes
r^*(\Omega^{\frac{1}{2}}\otimes E_0))$ and satisfy the relations of the 
lemma \ref{f3}. We denote these operators by $\nabla^{\E_0}$, $T^H$ and $l_0$.
The following linear spaces will be used later  
\begin{gather*}
\Omega^k\mathcal{E}_0=C_c^\infty(G^\partial,\Lambda^ks^*\nu^*\otimes
r^*(\Omega^{\frac{1}{2}}
\otimes(E_0\oplus E_0)))\quad;\quad\Omega^*\mathcal{E}_0=\oplus_k\Omega^k\E_0.
\end{gather*}
Here $\Omega^{\frac{1}{2}}$ refers to the longitudinal bundle of half densities 
of the foliation $(\partial M,\partial F)$, which is, in fact, the restriction of 
the longitudinal bundle $\bO$ to $(\partial M,\partial F)$, c.f. 
\cite[relation 4.47]{Me}. 
Now we define 
the $\Cl(1)$-superconnection $\mathbb B$, acting on $\Omega^*\mathcal{E}_0$, 
by the following formula 
(c.f, \cite[page 28]{MePi}).
\begin{gather}\label{csup}
\mathbb B=\alpha.\operatorname D_0+\Id.\nabla^{\E_0}-\frac{1}{4i}\alpha.\cl(T^H)
\end{gather}
where $\alpha=\begin{pmatrix}0&1\\1 &0\end{pmatrix}$ and 
$\Id=\begin{pmatrix}1&0\\0 &1\end{pmatrix}$ 
with respect to the obvious direct sum decomposition of $\Omega^*\mathcal{E}_0$. 
In fact $\mathbb B$ 
is the restriction of $\A$ to $\Omega^*\mathcal{E}_0$. 
Just like $\At$ in \ref{folman},  
the \Cl(1)-superconnection \Bb can be 
defined for the foliation 
$(\partial M,\partial F)$. In fact $\Bb$ is the restriction of 
the superconnection $\At$ to the boundary foliation $(\partial M,\partial F)$.
The invariance relation $\mathbb B^2-l=r^*{\Bb^2}$ holds and 
its proof is completely similar to the proof of the proposition 
\ref{li}. We
have $\mathbb B_s^2=s^2{\operatorname D}_0^2+sF_0+J_0$ where $F_0$ and $J_0$ are
differential operator with differential forms coefficient of positive degree.

\subsection{Algebraic structures of the holonomy groupoid}
Groupoid action of $G$ on itself provides an algebra structure on
\ob and a left \ob-module structure on \oE. Let $\phi$ and $\psi$ be
two elements of \ob and put
\begin{equation}\label{gg}
\phi.\psi\,(u):=\int_{G^{r(u)}}\phi(v)\wedge (v^*)^{-1}\psi(v^{-1}u)
\end{equation}
where $v^*$ denotes the dual of $v_*$. Notice that the integration is
 taken with respect to the density component of the first factor. Now
 let $\xi\in\oE$ and define
\begin{equation}\label{h}
\phi.\xi\,(u):=\int_{G^{s(u)}}\phi(v)(v^*)^{-1}\xi(uv).
\end{equation}
Here again the integration is taken with respect to the density
component of $\phi$.
With respect to these operations \oE is a left-module on the algebra
\ob. In what follows we characterize the structure of
$End_{\ob}(\oE)$ and describe its algebra structure. Let
$P:\E\rightarrow\oE$ be a linear operator which is defined through
its smooth kernel
 by 
\begin{gather}
P(\xi)(u)=\int_{G_{s(u)}}\xi(v)\,v^*K(v,u)~;\label{gol}\\
K(v,u)\in r^*_v(\wedge^*\nu^*\otimes\bO)\otimes r^*_u(\bO)\otimes
Hom(E_{r(v)},E_{r(u)}).\notag
\end{gather}
In above relation the integration is taken with respect to the product of 
half densities coming from $\xi$ and from the end point of $v$. 
It is easy to verify that $\B$-linearity is equivalent to the
following invariance property
\[K(v\gamma,u\gamma)=K(v,u)\quad\forall\gamma\in G^{s(u)}.\]
Therefore the operator $P$ may be defined as follows:
\begin{gather}
P(\xi)(u)=\int_{G_{s(u)}}\xi(v)\,v^*K(vu^{-1})~;\label{pso}\\
K(\gamma)\in r^*_\gamma(\Lambda^*\nu^*\otimes\bO)\otimes
s^*_\gamma(\bO)\otimes Hom(E_{r(\gamma)},E_{s(\gamma)}).\notag
\end{gather}
The operators defined in this way are called $G$-invariant. We denote the set
of all these operators by $Hom_G(\E,\oE)$.
By applying  $(vu^{-1})^*$ to $\xi(v)$ in relation \eqref{pso}, the
action of $P$ can be extended to $\oE$ and in this way $P$ defines an
element in $End_{\ob}(\oE)$. In particular we get the induced algebra
structure on $Hom_G(\E,\oE)$.

Notice that in relation \eqref{h} if $\xi$ belongs to 
$\Omega^*\mathcal{E}_0$ then $\phi.\xi$ belongs to $\Omega^*\mathcal{E}_0$ too.
So $\Omega^*\mathcal{E}_0$ is a right module on \ob. Let $Q$ be a 
linear operator on $\Omega^*\mathcal{E}_0$ defined through its smooth kernel
\begin{gather}
Q(\xi)(u)=\int_{G_{s(u)}^\partial}\xi(v)\,v^*K(v,u)~;\label{darti}\\
K(v,u)\in r^*_v(\Lambda^*\nu^*\otimes\Omega^{\frac{1}{2}})
\otimes r^*_u(\Omega^{\frac{1}{2}})\otimes
Hom(E_{r(v)},E_{r(u)})\notag.
\end{gather}
As in above, if the operator $Q$ is \ob-linear, i.e. 
$Q\in End_{\ob}(\Omega^*\mathcal{E}_0)$, then its kernel takes the 
form $K(v,u):=K_Q(vu^{-1})$. Since in this case $u$ and $v$ are in 
$G_{s(u)}^\partial$, we have $vu^{-1}\in G_\partial^\partial$,. 
So $Q$ may be defined through 
a kernel which is defined on $G_\partial^\partial$ 
\begin{gather}
Q(\xi)(u)=\int_{G_{s(u)}^\partial}\xi(v)\,v^*K_Q(vu^{-1})~;\label{dernh}\\
K_Q(\gamma)\in r^*(\Lambda^*\nu^*_\gamma\otimes\Omega^{\frac{1}{2}})\otimes
s^*_\gamma(\Omega^{\frac{1}{2}})\otimes Hom(E_{r(\gamma)},E_{s(\gamma)})~~~;
\gamma\in G_\partial^\partial.\label{dashy}
\end{gather}
Following \cite[relation 58]{GoLo}, we shall define a supertrace for 
operators $Q$ in $End_{\ob}(\Omega^*\mathcal{E}_0)$ with 
smooth kernels, i.e. the $G$-invariant smoothing operator on $\Omega^*\mathcal{E}_0$. 
Let $\tau$ be an element in $C_c^\infty(G,s^*(\text{}{^b\!}\Omega^1))$
satisfying the following relation (c.f. \cite[proposition 6.11]{Tu-l} 
for the proof of the
existence of $\tau$):
\begin{gather}
\int_{G^y}\tau(v)=1~;~~\forall y\in M.\label{tau0}
\end{gather}
Taking into the account the symmetries coming from the action of the
groupoid $G$ on $G^\partial$, following \cite[relation 58]{GoLo}, we define 
$STr(Q)\in C^\infty(M,\text{}{^b\!}\Omega^1\otimes\Lambda^*\nu^*)$ by  
\begin{gather}\label{tau}
STr(Q)(x)=\int_{G^\partial_x}\tau(v)\,\operatorname{str}(v^*K_Q(vv^{-1}))~.
\end{gather}
Here $\operatorname{str}$ is calculated with respect to the grading 
$r^*E_{|\partial M}=r^*E_0\oplus r^*E_0$ and 
the integration is taken with respect to the product of two 
half densities coming from the trace of the Schwartz kernel, c.f. 
relation \eqref{dashy}. Notice that  
in above relation $x$ is a point in $M$, so even though $Q$ is a 
boundary operator, $STr(Q)$ belongs to 
$C^\infty(M,\text{}{^b\!}\Omega^1\otimes\Lambda^*\nu^*)$. 
The following proposition is the counterpart of the proposition $2$ of 
\cite{GoLo} in our situation.
\begin{prop}\label{trace1}
\text{}\\
a) Let $\rho$ be a distributional linear functional on 
$C^\infty(M,\text{}{^b\!}\Omega^1\otimes\Lambda^k\nu)$
such that the linear functional $\eta$  defined on \ob  by
\[\eta(\phi)=\rho(\phi_{|M})\]
be a graded supertrace on \ob. Then $\rho\circ STr$ is independent of 
$\tau$ and defines a supertrace on 
$End_{\ob}(\Omega^*\E_0)$, i.e. vanishes on supercommutators.\\
b)If in addition  $\eta$ is a closed graded trace on $\Omega^k\B$,
i.e. $\eta\circ d^H=0$, then for all 
$Q\in End_{\ob}(\Omega^*\E_0)$ 
\[\rho\circ STr[\nabla^{\E_0},Q]=0\]
\end{prop}
\begin{pf}
It suffices to prove this proposition for $Q\in End_\B(\E_0)$, 
the general case follows easily 
from this special case. 
We show at first the independence from $\tau$. 
Let $\chi(G_\partial)$ denote the characteristic function of
$G_\partial$ which can be approximated by continuous 
functions. From \eqref{gg} we have 
\[\chi(G_\partial).\phi(x)
=\int_{G^x}\chi(G_\partial)(v)\wedge\phi(v^{-1})
=\int_{G_\partial^x}\phi(v^{-1}).\]
On the other hand it is easy to verify that $\phi.\chi(G_\partial)(x)=0$ if 
$x\notin\partial M$ while 
\[\phi.\chi(G_\partial)(x)=\int_{G^x}\phi(v)\chi(G_\partial)(v^{-1})
=\int_{G^x}\phi(v)~;\text{ for }x\in \partial M.\]
The trace property of  $\eta$ reads 
$\eta(\chi(G_\partial).\phi)=\eta(\phi.\chi(G_\partial))$ for each 
$\phi\in \B$, so the above relations imply the following one
\begin{equation}\label{char}
\int_M\rho(x)\int_{G_x^\partial}\phi(v)=\int_{\partial
M}\rho(y)\int_{G^y}\phi(v).
\end{equation}
Here we have used the fact that $\rho$ is distributional 
and can be approximated by distributions defined by 
smooth functions $\rho(x)$ (e.g. see \cite[theorem 6.3]{LiebLoss-analysis}).
If we replace $\chi(G_\partial)$ by $1-\chi(t\geq \epsilon)$ we
obtain
\begin{equation}\label{char1}
\int_M\rho(x)\int_{G_x^{t\geq\epsilon}}\phi(v)
=\int_{M^{t\geq\epsilon}}\rho(y)\int_{G^y}\phi(v).
\end{equation}
Here $t\geq0$ is a $G$-invariant defining function for the boundary
of leaves $G_x$'s. Now let $Q$ be a smoothing $G$-invariant operator 
with kernel $K_Q$. 
Using \eqref{char} we obtain the following equalities 
where $O_y^y$ denotes the unite elements of 
the group $G_y^y$ for each $y\in M$
\begin{align}
\rho\circ STr\, Q&=\int_M\int_{G_x^\partial}\rho(x)\tau(v)\,strK_Q(vv^{-1})\notag\\
&=\int_{\partial M}\rho(y)\int_{ G^y}\tau(v)\,strK_Q(vv^{-1})\notag\\
&=\int_{\partial M}\rho(y)\,str\,K_Q(O_y^y)
\int_{G^y}\tau(v)\notag\\
&=\int_{\partial M}\rho(y)\,str\,K_Q(O_y^y).\label{foot}
\end{align}
The last expression does not depend on $\tau$ therefore $ \rho\circ STr\,Q$
is independent of $\tau$. Moreover this formula can be 
used to show that $\rho\circ STr$ is a supertrace. For this purpose it suffices 
to show that $\rho\circ STr[Q,Q']=0$ when $Q,~Q'$ belong to $End_\B(\E_0)$ with 
$\dim_\C E_0=1$. In this situation $K_Q$ and $K_{Q'}$ are compactly supported 
smooth functions 
on $G_\partial^\partial$ and we consider them as $G$-invariant compactly 
supported functions on $G$, i.e. elements in $\B$. By \eqref{gg} and \eqref{dernh} 
we have 
\[K_{QQ'}(O_y^y)=\int_{G_\partial^y}K_{Q'}(\gamma)K_Q(\gamma^{-1})=(K_{Q'}.K_Q)(O_y^y)\]
where $K_{Q'}.K_Q$ denotes the multiplication of two elements of the algebra $\B$. 
Using \eqref{foot} and above relation we obtain  
\begin{align*}
\rho\circ STr[Q,Q']&=\int_{\partial M}\rho(y)(K_{Q'}.K_Q-K_Q.K_{Q'})(O_y^y)\\
&=\eta((K_{Q'}.K_Q-K_Q.K_{Q'})_{|M})
\end{align*}
But the last expression vanishes by the assumption on $\eta$ and this 
completes the proof of the first part of the proposition.

Part (b) of this proposition follows from the following relation that we 
are going to prove.
\begin{equation}\label{noh}
\rho\circ STr[\nabla^{\E_0}, Q]=\eta(d^H\circ r^*(STr\,Q)).
\end{equation}
Using the Leibnitz rule for $\nabla^{\E_0}$, it suffices to prove
this relation for $Q\in End_G(\E_0)$.
Let $T_1^HG$ be a connection which is integrable in an small neighborhood $V$
of $v$ (that is $T_1^G$ is the tangent bundle of local 
transversal sub manifolds of $M$). 
A horizontal vector $\mathbf y_1\in T_{1v}^HG$ determines uniquely a 
horizontal vector $\mathbf y\in T_v^HG$ such that $\mathbf y_1-\mathbf y$ 
is a vertical vector. 
Denote by $d_1^H$ and $\nabla^{\E_0}_1$ the associated operators to the 
integrable connection $T_1^HG$. It follows from the above discussion that 
$d_1^H-d^H$ is the differentiatian with respect to vertical vectors which implies 
the vanishing of $(d_1^H-d^H)r^*(Tr\,Q)$.
On the other hand,
$\nabla^{\E_0}_1-\nabla^{\E_0}$ is an element of  
$Hom_G(\E_0,\Omega^*\E_0)$ which implies  
$STr[\nabla^{\E_0}_1,Q]=STr[\nabla^{\E_0},Q]$.
Thus to prove relation \eqref{noh}, we may assume that $T^HG$ is integrable.
 So one can find local coordinate system $T\times V$ around $v\in G$ such that 
$T^hM$ is the tangent bundle of transversal $T$. Denote the coordinate of $T$ 
by $z$. The horizontal differentiation is nothing else than the 
partial derivatives $\partial_z$. 
Considering $K_Q$ as a local matrix-valued function it is easy to see that 
$str[\partial_z,K_Q](v)=\partial_z(str\,K_Q(v,v))$, that means 
$str[\nabla^{\E_0},K_Q](v)=d^H(str\,K_Q)(v)\in\Lambda^*\nu^*$ which with 
\eqref{pulinv} and \eqref{char} imply the following equalities and complete the proof
\begin{align*}
\rho\circ STr[\nabla^{\E_0}, K_Q]
&=\int_M\rho(x)\int_{G_x^\partial}\tau(v)\,str\,[\nabla^{\E_0},Q](vv^{-1})\\
&=\int_M\rho(x)\int_{G_x^\partial}\tau(v)\,d^H(str\,K_Q)\,(vv^{-1})\\
&=\int_{\partial M}\rho(y)\,d^h(str\,K_Q(O_y^y))\int_{G^y}\tau(v)\\
&=\int_{\partial M}\rho(y)\,d^h(\int_{G^y}\tau(v)\,str\,K_Q(vv^{-1}))\\
&=\int_M\rho(x)\,d^h(\int_{G_x^\partial}\tau(v)\,str\,K_Q(vv^{-1}))\\
&=\eta(d^H\circ r^*(STr\,Q))
\end{align*}
\end{pf}
\begin{ex}
Let $\mathcal{C}^k$ be a holonomy invariant $k$-current acting on
$C^\infty(M,\Lambda^k\nu^*)$ and let $\{U_i\}_i$ denote a finite
covering for $(M,F)$ by the flowboxes $U_i=T_i\times V_i$ where
$\dim(V_i)=p$. 
Denote by $\mathcal{C}^k_i$ the restriction of $\mathcal{C}^k$ to $T_i$. Let
$\{\phi\}_i$ be  a partition of unity subordinate to this covering.
The linear functional $\rho$ given by   
\begin{gather}
\rho:C^\infty(M,\text{}{^b\!}\Omega^1(F)\otimes\Lambda^*\nu^*)\rightarrow \C\notag\\
\rho(\omega)=\sum_i <\int_{V_i}\phi_i\omega,\mathcal{C}^k_i>\label{majid}
\end{gather}
satisfies the condition of the above proposition, c.f. 
examples 6 and 7 of \cite{GoLo}. In above the integration on $V_i$  
is taken with respect to the $\text{}{^b\!}\Omega^1$-factor 
of $\omega$ while the integration on $V_i\cap\partial M$ is 
taken with respect to the $\text{}{^b\!}\Omega^1$-factor 
of $\omega$.
\end{ex}

Now let $Q$ be as in above acting on   
$C_c^\infty(G^\partial,\Lambda^ks^*\nu^*\otimes
r^*(\Omega^{\frac{1}{2}} \otimes(E_0\oplus E_0)))$ with smooth kernel 
$K$ such that 
\begin{equation}\label{gali}
K(v,v)=\Id\,K_0(v,v)+\alpha\,K_1(v,v)
=\left(\begin{array}{cc}K_0&K_1\\K_1&K_0\end{array}\right)
\end{equation}
where $K_0(v,v)$ and $K_1(v,v)$ are elements in 
$\Lambda^*\nu_x^*\otimes\Omega^1(\text{}{^b\!}T_yF)\otimes End(E_0)$ 
with $x=s(v)$ and $y=r(v)$, while   
$\Id=\begin{pmatrix}1&0\\0 &1\end{pmatrix}$ and 
$\alpha=\begin{pmatrix}0&1\\1 &0\end{pmatrix}$. 
Put $str_\alpha K(v,v):=tr(\Id\,K_1)$ and define 
the $\CCl(1)$-supertrace(or $\alpha$-supertrace) of $Q$, denoted by 
$STr_\alpha(Q)\in C^\infty(M,\text{}{^b\!}\Omega^1\otimes\Lambda^*\nu^*)$, 
by the following formula 
\begin{gather}\label{clsuptr}
STr_\alpha(Q)(x):=\int_{G^\partial_x}\tau(v)\,\operatorname{str_\alpha}
(v^*K_{|G_x^\partial}(v,v))~.
\end{gather}
It is clear from above discussion that if $Q$ is an even operator on 
$C_c^\infty(G^\partial,\Lambda^ks^*\nu^*\otimes
r^*(\Omega^{\frac{1}{2}} \otimes(E_0\oplus E_0)))$
then $STr_\alpha(Q)$ belong to 
$C^\infty(M,\text{}{^b\!}\Omega^1\otimes\Lambda^{odd}\nu^*)$ 
\begin{rem}
 Let $r^*(\omega)$ be a $G$-invariant element of $\ob$. The proof of
 the relation (\ref{foot}) can be applied to prove the following
 relation which will be used later
\begin{equation}\label{evaa}
\int_M\int_{G_x}\rho(x)\tau(v)
r^*(\omega)(v)=\int_M\rho(x)\omega(x).
\end{equation}
\end{rem}
Similar to the formula \eqref{tau} we could define $STr(P)$ for 
$P\in End_\B(\E)$ with compact support disjoint from boundary. This condition 
on the support cannot be satisfied always, so we have to define a more 
general trace for a larger
class of $G$-invariant operators. These are smoothing $G$-invariant
b-pseudodifferential operators which are introduced by R. Melrose \cite{Me}. 
In the following section we follow \cite{MePi} and \cite{LePi-b}
to describe those aspects of Melrose's b-calculus which are
necessary for future uses in this work.

\section{Some aspects of G-invariant $b$-calculus}\label{section3}
The underlying space for small b-calculus is an appropriate
b-stretched product. For each $x\in M$, the leaf $F$, passing through $x$, 
with its exact b-metric is a b-manifold. The covering space with boundary
$G_x\stackrel{r}{\rightarrow} F_x$, equipped with the pull-back metric is 
a non compact b-manifold. Let $B_x$ consists of those elements $(u,v)$ of 
$G_x^\partial\times G_x^\partial$ such that $u$ and $v$ belong to the 
same connected component of $G_x^\partial$. The $b$-stretched product 
$G_{x,b}^2$ is obtained from $(G_x)^2$ by replacing $B_x$ with 
inward-pointing unit normal bundle to $B_x$, i.e. $S_+N(B_x)$, and putting 
an appropriate topology on it. $S_+N(B_x)$ is denoted by $\ff_x$ and is 
called the front face of b-stretched product space $G_{x,b}^2$. 
It turns out that this space is equipped with a surjective smooth blow-down map 
$(\beta_x)_b:G_{x,b}^2\to (G_x)^2$ which restricts to a diffeomorphism 
$G_{x,b}^2\verb+\+S_+N(B_x)\to(G_x)^2\verb+\+B_x$. Using this diffeomorphism 
one can lift the subset $G_x^\partial\times \dot G_x$, $G_x\times G_x^\partial$ 
and $\triangle_x=diag(G_x\times G_x)$ of $(G_x)^2\verb+\+B_x$ to $G_{x,b}^2$. 
We denote the closure of the lifted sets, respectively, by $\lb_x$, $\rr_x$ 
and $\triangle_{x,b}$. These are left boundary, right boundary and 
the b-diagonal of b-stretched product space $G_{x,b}^2$. The space $G_{x,b}^2$ 
is the carrier for the kernel of b-pseudodifferential operators on $G_x$.
 
As in the family case (cf. the appendix of \cite{MePi}), the
b-stretched product of the holonomy groupoid $G$ is the smooth family, 
parameterized by $x\in M$,
of stretched product spaces $G^2_{x,b}$. We denote this space by
$G_b^2$. Clearly each element of $G_x^y$ provides a diffeomorphism from 
$G_{y,b}^2$ to $G_{x,b}^2$. An object defined on $G_b^2$ is called 
$G$-invariant if it is invariant with respect to this action of $G$ 
on $G_b^2$. We denote the family $\beta_{x,b}: G^2_{x,b}\rightarrow G_x^2$ of
blow-down maps by $\beta_b: G_b^2\rightarrow G^2$. Similarly we
denote by $\lb$, $\rr$ and $\ff$ respectively
the $G$-invariant left, right, and the front face family of $G_{x,b}^2$'s, while
$\triangle_b$ denotes the $G$-invariant family of b-diagonals $\triangle_{x,g}$. 
The function $r:=x+x'$ is a $G$-invariant defining function for $\ff$. 
For each $x\in M$, the small b-calculus $\Psi^*_b(G_{x,b}^2,r^*(\bO\otimes E))$ 
is defined exactly 
as in \cite[chapter 4]{Me} through a precise description of the kernels 
on $G_{x,b}^2$. By a $G$-invariant b-pseudodifferential operator $P$, we mean 
a family $P_x$ of b-pseudodifferential operators with kernels $K_x$ such 
that $K_x$ is $G$-invariant. $P$ is of order $m$ if for each $x$ the 
operator $P_x$ is of order $m$. Similarly $P$ is classical or one-step 
polyhomogeneous if each $P_x$ is one step polyhomogeneous. We denote the set 
of all $G$-invariant b-pseudodifferential operators by 
$\Psi^*_{b,G}(G_b^2,r^*(\bO\otimes E))$ and the set of all classical 
$b$-pseudodifferential operator by 
$\Psi^*_{b,G,os}(G_b^2,r^*(\bO\otimes E))$. The operator 
$P$ is elliptic if $P_x$ is elliptic for each $x\in M$. 
For each $x\in M$ the manifold $G_x^\partial$ 
is smooth and without boundary 
and $\mathcal E_0=r^*(E_{|\partial M})$ is a vector bundle over it. The 
diagonal action of each element of $G_x^y$ provides a diffeomorphism from 
$G_y^\partial\times G_y^\partial$ to $G_x^\partial\times G_x^\partial$. 
Now we can proceed as in above and define the set of all 
$G$-invariant pseudodifferential operators 
$\Psi^*_{G}(G^\partial,r^*(\Omega^{\frac{1}{2}})\otimes \mathcal E_0)$.
Unlike to the compact case, the small
b-calculus $\Psi^*_{b,G}(G_b^2,\bO)$ and the pseudodifferential calculus 
$\Psi^*_{G}(G^\partial,r^*(\Omega^{\frac{1}{2}})\otimes \mathcal E_0)$ 
are not algebras. Consequently
the principal symbol exact sequence can not be used to construct a
parametrix for, e.g. an elliptic b-differential operator. To overcome this
difficulty, following \cite[definition 4.6]{LePi-b}, we make some assumption on the
behavior of kernels far from b-diagonal, and also on the nature of
the foliation groupoid.

Let $0<\epsilon<1$ and set
\[\mathcal{O}_\epsilon(\ff)=\cup_{x\in M}
\{u\in G_{b,x}^2|d(\beta_{x,b}(u),B_x)<\epsilon\}\]
Here and in what follows  $d(,)$ denotes the $G$-invariant distance
function coming from an ordinary $G$-invariant metric $\tilde g$ on $G$. 
Let $x$ and $x'$ be two $G$-invariant defining functions, respectively, 
for $\lb$ and $\rr$ and put $r:=x+x'$ and $\tau:=(x-x')/(x+x')$ 
If $(v,v')$ is a local coordinate for
$B:=\cup_xB_x$ then $(r, \tau, v, v')$ define a
local coordinate system for $\mathcal{O}_\epsilon(\ff)$ around $\ff$.

\begin{defn}\label{rapdec}
Let $K\in C^{-\infty}(G_b^2,\bO)$ be a distributional section of
$r^*\bO\otimes r^*\bO$ which is smooth outside a small neighborhood of 
b-diagonal $\triangle_b$. One says $K$ has the
rapidly decreasing property outside $\epsilon$-neighborhood of 
lifted diagonal $\triangle_b$ if \\
i) For each multi-index of derivations $\alpha$ and any $N\in\N$
there is a constant $C_{\alpha,N}$ such that for all $(u,u')\in
G_b^2\verb+\+\mathcal{O}_\epsilon(\ff)$ satisfying
$d(u,u')>\epsilon$ one has
\begin{equation}\label{fc}
|\nabla^\alpha K(u,u')|(1+d(u,u'))^N<C_{\alpha,N}
\end{equation}
ii) There exists a constant $D_{\alpha,N}$
such that for each $(r,\tau,v,v')$ with $d(v,v')>\epsilon$ one
has
\[|\nabla^\alpha K(r,\tau,v,v')|(1+d(v,v'))^N<D_{\alpha,N}\]
\end{defn}
Notice that, due to the compactness of $M$, this definition is
independent of the $G$-invariant metric.
For each $m\in \R\cup\{-\infty\}$ denote by $\RR^m(G;r^*E)$ the subspace 
of  $\Psi^m_{b,G}(G;r^*(\bO\otimes E))$
consisting of the operators with rapidly decreasing kernels outside 
of a $\epsilon$-neighborhood of $\triangle_b$, and put
$\RR^*(G,r^*E)=\cup_m\RR^m(G;r^*E)$. 
Similarly define $\mathcal R_G^m(G^\partial;\mathcal E_0)$ 
as the subset of  
$\Psi^m_{G}(G^\partial,\mathcal E_0)$ consisting of 
operators with rapidly 
decreasing kernels outside a $\epsilon$-neighborhood of the diagonal 
of $G^\partial\times G^\partial$ and put 
\[\mathcal R_G^*(G^\partial;\mathcal E_0)
=\cup_m\mathcal R_G^m(G^\partial;\mathcal E_0)~;\hspace{1cm}
\mathcal R_G^{-\infty}(G^\partial;\mathcal E_0)
=\cap_m\mathcal R_G^m(G^\partial;\mathcal E_0)~.\]
The spaces $\RR^*$ and $\mathcal R_G^*$ are algebras if 
for each $x\in M$ and $N\in\N$  sufficiently large, the following 
integrals are convergent  
\[\int_{G_x\times G_x}\theta(y',y)(1+d(y',y))^{-N}\,d\mu_{\hat{g}}
~,\hspace{5mm}
\int_{G_x^\partial\times G_x^\partial}\theta'(y',y)(1+d(y',y))^{-N}\,d\mu_{\hat{g}}.\]
Here $\theta$ and $\theta'$ are respectively the characteristic function of the 
subsets of $G_x\times G_x$ and $G_x^\partial\times G_x^\partial$  consisting of 
the points  $(y',y)$ with $d(y',y)\ge\epsilon$ for a positive real number $\epsilon$. 
This is the case if the Riemannian manifolds $(G_x,\hat{g})$ 
and $(G_x^\partial,\hat{g})$ 
are of polynomial growth (See \cite[section 2]{Gr-g}).
Let $\hat{X}$ be a Riemannian manifold (non compact and
with boundary) and let $\Gamma$ be a discrete group of 
isometries of this manifold which acts properly on $\hat{X}$ such
that $X:=\hat{X}/\Gamma$ is a Riemannian manifold. The proof of
\cite[lemma 2]{Mi-a} shows that if $X$ and $\Gamma$ are of
polynomial growth then $\hat{X}$ is of polynomial growth too. This
discussion shows that under the following hypothesis $G_x$'s are of
polynomial growth, so the spaces $\RR^*(G,r^E)$ and
$\mathcal{R}_G^*(G^\partial,r^*E_{|\partial M})$ have
the structures of algebras. In particular the spaces of  rapidly decreasing
smoothing operators $\mathcal{R}^{-\infty}_{b,G}(G,r^*E)$ 
and $\mathcal{R}_G^{-\infty}(G^\partial,r^*E_{|\partial M})$ have the 
algebra structure .
>From now on the following hypothesis is assumed to
be satisfied.
\begin{hyp}\label{hyp2}\text{}\\
i) The leaves of the foliations $(M,F)$ and $(\partial M,\partial F)$ are 
of polynomial growth with
respect to a (hence to all) smooth longitudinal Riemannian metric on $M$ 
and $\partial M$.\\
ii) For each $x\in M$ the holonomy group $G_x^x$ is of polynomial
growth. Moreover, the actions of $G_x^x$ on $G_x$ and $G_x^\partial$ are 
proper.
\end{hyp}
\begin{lem}\label{frd}
The small space of $G$-invariant rapidly decreasing
classical b-pseudodifferential operators $\mathcal{R}^*_{b,G}(G;r^*E)$ is an
algebra. Moreover, one has the following exact sequence where
$\text{}{^b\!}\sigma_m$ denotes the principal symbol map.
\begin{equation}\label{dell}
0\rightarrow \mathcal{R}^{m-1}_{b,G,os}\rightarrow
\mathcal{R}^m_{b,G,os}\stackrel{\text{}{^b\!}\sigma_m}\rightarrow
C_G^\infty(S(\text{}{^b\!}T^*G) ; End(\pi_b^*(r^*E)))\rightarrow0.
\end{equation}
\end{lem}
\begin{pf}
This short sequence is exact when is restricted to $G_x$ for any $x\in M$. 
In despite of the non compactness of the underlying manifold $G_x$, this 
follows from \cite[proposition 4.23]{Me}. 
\end{pf}
 As in the ordinary pseudodifferential calculus, the above lemma implies the
existence of parametrices for elliptic b-operators. More precisely,
given an elliptic operator $P\in \mathcal{R}^m_{b,G,os}(G,r^*E)$,
there exists $Q\in \mathcal{R}^{-m}_{b,G,os}(G,r^*E) $ and the
smoothing operator $R$ in $\mathcal{R}^{-\infty}_{b,G}(G,r^*E)$ such
that
\begin{equation}\label{ain}
P\circ Q-Id=R.
\end{equation}
Moreover, since the principal symbol depends only on the behavior of kernel 
near $\triangle_b$, the kernel of the parametrix $Q$ may be assumed 
supported in a small neighborhood of $\triangle_b$. 
If $P$ is a differential operator this implies that the kernel of $R$ 
is also supported in a small neighborhood of $\triangle_b$. 
Let $P$ be a $G$-invariant operator  with rapidly decreasing  
kernel $K_P$ in the small b-calculus, $P\in\mathcal R^*_{b,G}(G,r^*E)$. 
The boundary operator of $P$ is the operator 
$P_{|G^\partial}\in\mathcal R^*_{G}(G^\partial,r^*E_{|\partial M})$ 
defined by the following relation
\begin{equation}\label{frt}
P_{|G^\partial}(\xi_0)=P(\xi)_{|G^\partial}.
\end{equation}
Here $\xi_0$ is an element 
in $C_c^\infty(G^\partial,r^*E_{|\partial M})$ and 
$\xi$ is an arbitrary extension of $\xi_0$ to an element in 
$C_c^\infty(G,r^*(\bO\otimes r^*E))$. Since the kernel of operators in 
small b-calculus vanish smoothly at left and right boundaries of $G^2_b$,
the above definition is independent of the extension $\xi$. This implies also
the product formula 
$(P\circ Q)_{G^\partial}=P_{|G^\partial}\circ Q_{|G^\partial}$. 
We recall that $r=x+x'$ is a defining function for the front face $\ff$ in 
$G_b^2$. If the boundary operator of $P$ vanishes then the kernel $K_p$ of 
$P$ has to vanish on $\ff$. So $K_P$ takes the form $rK_{P'}$ where 
$P'\in \mathcal R^*_{b,G}(G,r^*E)$. This follows from the explicite 
description of the behavior of kernels near $\ff$, see 
\cite[relation 4.61]{Me}. 

>From now on we fix a G-invariant co-normal structure 
$n$ for $G^\partial\subset G$. This consists of a $G$-invariant family 
$\{n_x\}_{x\in M}$ of sections of    
the normal bundle of $G_x^\partial$ in $G_x$.  
Let $\sigma$ be a defining function for $G^\partial\subset G$ 
satisfying $d\sigma(n)=1$.
Following \cite[Proposition 5.8]{Me}, we define the indicial
family $I_n(P,\lambda)\,;\lambda\in\C$ of a G-invariant
b-pseudodifferential operator $P$ by the following formula which is 
independent of $\sigma$ as long as it satisfies the condition $d\sigma(n)=1$
\begin{equation*}
I_n(P,\lambda)=(\sigma^{-i\lambda}P\,\sigma^{i\lambda})_{|G^\partial}.
\end{equation*}
We have clearly the product formula $I_n(P\circ
P',\lambda)=I_n(P,\lambda)\circ I_n(P',\lambda)$. As an example one has 
\begin{equation}\label{ind12}
I_n(\A,\lambda)=\lambda.\gamma+\mathbb B~;\hspace{1cm}
\gamma:=\left( \begin{array}{cc}
0&-i\\i&0
\end{array}\right).
\end{equation}
It is clear from above discussion that 
$I_n(P,\lambda)=0$ if and only if the kernel $K_p$ of $P$ vanishes on $\ff$.\\
We recall from \eqref{ain} that each $G$-invariant
b-differential operator $P$ has a $G$-invariant parametrix $Q$ in
small b-calculus such that the remainder $R=P\circ Q-Id$ is a smoothing
b-pseudodifferential operator. In what follows we need to have
a good parametrix with remainder having a vanishing indicial operator. Such a
parametrix does not exist in small b-calculus but in a slightly
extended calculus which is, the so-called, full calculus with bounds. 
We recall from \cite[section 7]{LePi-b} that an operator $P$ in the 
full calculus with bound $\mathcal{R}^{*,\delta}_{b}(G_x,r^*(\bO\otimes E))$ 
is defined by a precise description of the singularities of its kernel 
on $G^2_x$. In this extended calculus, the singularity on $\triangle_b$ 
is the same as in the small b-calculus. This part of singularities 
determines the symbol and the degree of the operator.   
In contrast to the small b-calculus, $K_P$ is singular near the left, the right 
and the front face hypersurfaces of $G_{x,b}^2$. 
The number $\delta$.which determines the order of these singularities 
is determined by $I_n(P,\lambda)$. Let $\delta$ be a non negative number and
the operators satisfy the rapidly decreasing property of 
definition ref{rapdec}. Then, with the hypothesis 
\ref{hyp2}, the full calculus 
$\mathcal{R}^{*,\delta}_{b}(G_x,r^*(\bO\otimes E))$ is an algebra.
The full calculus with bounds $\mathcal{R}^{m,\delta}_{b,G}$ consists of 
the $G$-invariant family of full calculus for $G_x$'s. 
To show the existence of a parametrix with remainder having a vanishing indicial
operator we will suppose the following hypothesis
\begin{hyp}\label{hyp}
The family of Dirac operators $D_0$ is $L^2$-invertible with
bounded inverse. In other word, there exists $\epsilon_0>0$ such that
for each $x\in M$, the $L^2$-spectrum of the unbounded operator 
 \[D_{0|G_x}:L^2(G_x^\partial;
r^*E_{|\partial M})\to L^2(G_x^\partial;r^*E_{|\partial M})\]
is disjoint from $[-\epsilon_0,\epsilon_0]$.
\end{hyp}
Notice that this is a global hypothesis depending on the geometry of
$G^\partial$. The indicial family of $D$ is given 
by $I_n(D,\lambda)=\lambda\gamma+\alpha D_0$.
So, as a consequence of the above hypothesis, $I_n(D,\lambda)$ is invertible if 
$-\epsilon_0\leq\rm{Im}(\lambda)\leq\epsilon_0$.

\begin{prop}\label{frad}
Let the above hypothesis be satisfied by the boundary operator
$D_0$ and let $0<\delta<\epsilon_0$. Then
there exists a parametrix $Q\in \mathcal{R}^{-1,\delta}_{b,os,G}$ such
that
\[DQ-Id=R~;\hspace{5mm}QD-Id=R',\]
where $R$ and $R'$ are in $\rho_{\ff}\mathcal{R}^{-\infty,\delta}_{b,G}$. 
In particular $I_n(R,\lambda)=I_n(R',\lambda)=0~$. 
\end{prop}
\begin{pf}
Let $T$ be a small transversal sub manifolds of $M$.
We refer to \cite[theorem 7.1]{LePi-b} for the construction of the 
smooth family of parametrices $Q_x$ for $D_{|G_x}$, $x\in T$ such that 
\[D_xQ_x-Id=R_x~;\hspace{5mm} 
R\in \rho_{\ff}\mathcal{R}^{-\infty,\delta}_{b}(G_x,r^*(\bO\otimes E)).\]
In addition to the properties described in 
the proposition, the kernel of $Q_x$ and of $R$ are invariant with respect 
to the action of $G_x^x$ on $G_x$. Let $y\notin T$ and 
take $u\in G_x^y$ which provides a diffeomorphism $u:G_x\to G_y$. This 
diffeomorphism can can be used to 
define  $Q_y:=u^*Q_x$ and $R_y:=u^*R_x$ which belong, respectively, to 
$\mathcal{R}^{-1,\delta}_{b}(G_y,r^*(\bO\otimes E))$ 
and $\rho_{\ff}\mathcal{R}^{-\infty,\delta}_{b}(G_y,r^*(\bO\otimes E))$ 
such that $D_yQ_y-Id=R_y$. 
In this way we have extended the family of parametrices $Q_x$, $x\in T$ 
to a family on the set $G_T\cup G^T$. The groupoid $G$ can be covered 
by disjoint union of such sets. We can apply the above method to 
define, $G$-invariantly, the parametrix family on whole $M$.
\end{pf}
For each $x\in M$, let $\D_x$ denote the restriction $\D_{|G_x}$.
For $s>0$, one can define the heat operator $e^{-s^2\D_x}$
which belongs to $\mathcal R^{-\infty}_{b,G}(G_x; r^*(\bO\otimes E))$. 
To see the construction of this heat operator 
we refer to 
\cite[Theorem 10.3]{LePi-b}. Putting together all these operators for all 
$x\in M$ we obtain the $G$-invariant heat operator $e^{-s^2\D}$. 
This operator is an element of 
$\mathcal R^{-\infty}_{b,G}(G; r^*(\bO\otimes E))$ provided that $s>0$. 
Now, using the Volterra
formula (cf. \eqref{ik0} and \eqref{ik}) we can define for $s>0$ the
heat kernel $e^{-(\As^2-l)}$ of the Bismut superconnection $\As$.
This operator belongs to 
$\mathcal R^{-\infty}_{b,G}(G;\,r^*(\bO\otimes\Lambda^*\nu^*\otimes E))$.

Let $\pi_F\colon\text{}{^b\!}TM\to \text{}{^b\!}TF$ be the projection on 
the longitudinal b-tangent bundle with kernel $T^hM$. 
Using a scalar product $g_h$ on $T^hM$ one get the direct sum scalar product 
$g_h\oplus g$ on the bundle $\text{}{^b\!}TM$. Let $\tilde{\nabla}$ be
the Levi-Civita connection associated to this scalar product and put
\[\tilde{\nabla}^{M/F}=\pi_F\tilde{\nabla}\,\pi_F.\]
The connection $\tilde{\nabla}^{M/F}$ for $\text{}{^b\!}TF$ is  
independent of the horizontal metric $g_h$, c.f. \cite[chapter 10]{BeGeVe}. 
We denote by
$\tilde{R}^{M/F}$ the curvature of this connection which is a smooth section of 
$End(TF))$ with coefficients in $\Omega^2(M)$. The
longitudinal $\hat{A}$-genus of the foliation $(M,F)$ is an element of 
$\Omega^*(M)$ defined by
\begin{equation}\label{akar}
\hat{A}(M,F)=\operatorname{det}\left(\frac{\tilde{R}^{M/F}/2}
{\sinh(\tilde{R}^{M/F}/2)}\right)^{1/2} 
\end{equation}
Before stating the next proposition, we introduce some notation and concepts.
We recall that $E$ is a longitudinal Clifford bundle. 
It turns out that $E$ has a Clifford connection $\nabla^{E/S}$ which is 
compatible with the Levi-Civita connection
$\tilde{\nabla}^{M/F}$. Using Proposition 3.43 of \cite{BeGeVe} the
curvature of this connection $\tilde R^{E/S}$ has the following decomposition
\begin{equation*}
\tilde R^{E/S}=\operatorname{cl}(\tilde R^{M/F})+\tilde F^{E/S}
\end{equation*}
where $\operatorname{cl}(\tilde R^{M/F})$ is the action of the Riemannian
curvature of the connection $\tilde{\nabla}^{M/F}$ on the bundle $E$
through the Clifford action. More precisely we have
 \[\operatorname{cl}(\tilde R^{M/F})(e_\alpha,e_\beta)=\frac{1}{4}
\sum_{1\leq i,j\leq p}(\tilde{R}^{M/F}(e_\alpha,e_\beta)e_i,e_j)
\operatorname{cl}(e_i)\operatorname{cl}(e_j)\]
where $\{e_\alpha\}_\alpha$ is a basis for the
tangent bundle $TM$ while $\{e_i\}_i$ is a basis for the
longitudinal tangent bundle $TF$. The  term $\tilde F^{E/S}$ is 
an invariant of the Clifford connection
$\nabla^{E/S}$ which is called the  twisting
  curvature of $\nabla^{E/S}$.
The twisted Chern character of the vector bundle $E$ is a differential form 
on $M$ defined by
\begin{equation}\label{twi}
\operatorname{Ch}(E/S)=\operatorname{str\,e}^{-\tilde F^{E/S}}.
\end{equation}
If $E=S\otimes V$ then the twisting curvature is the curvature of the 
twisting hermitian bundle $V$. Notice that because at the boundary 
point  $\tilde\nabla_{\vec n}=0$, the differential forms $\hat{A}(M,F)$ 
and $\operatorname{Ch}(E/S)$ are tangential at boundary points. 
Using the decomposition (\ref{ff}) one defines the following projection
\[\{~.~\}_p:\Lambda^*(\text{}{^b\!}TM)
\rightarrow\Lambda^*(T^{h*}M)\otimes\Lambda^p(\text{}{^b\!}TF)
\stackrel{\simeq}{\rightarrow}\Lambda^*(T^{h*}M)\otimes\text{}{^b\!}\Omega^1(F).\]
In above the isomorphism between the longitudinal bundle 
$\Lambda^p(\text{}{^b\!}TF)$ and the bundle of longitudinal $b$-densities of 
order $1$ is given by the longitudinal $b$-metric $g$.
\begin{prop}[$G$-invariant Bismut b-density theorem]\label{lim}
We have the following asymptotic relation for the restriction of the heat
 kernel of the Bismut superconnection to the 
diagonal $\triangle_b$. This relation occurs in
 $\Lambda^*(T_u^HG)^*\otimes\Omega^1(F)$
\begin{equation}\label{limmm}
\lim_{s\rightarrow 0}\operatorname{str}\operatorname{e}^
{-(\As^2-l)}(u,u)
=\frac{1}{(2\pi i)^{p/2}}r^*\{\operatorname{\hat{A}}(M,F)
\operatorname{Ch}(E/S)(y)\}_p
\end{equation}
where $r(u)=y$.
\end{prop}
\begin{pf}
 Proposition \ref{li} asserts that $\As^2-l=r^*\Ats$, so 
$e^{-(\As^2-l)}=r^*(e^{\Ats})$.
This proves the $G$-invariance of the restriction of the 
kernel $e^{-(\As^2-l)}$ to the diagonal $\triangle_b(G)$. So
\[\operatorname{str}\,e^{-(\As^2-l)}(u,u)=r^*(\operatorname{str}
\,e^{-\Ats^2}(y,y))\]
where $y=r(u)$. Since $(M,F)$ is locally a fibration of 
b-manifolds and the assertion 
of the proposition is local, we may assume that $(M,F)$ is a fibration 
of $b$-manifolds and $E=S\otimes V$. For a fixed $y\in F$, 
the heat kernel $k_s(y,y):=e^{-\Ats^2}(y,y)$ 
is a smooth $s$-depending element in 
$\Omega^1_y\otimes \Lambda^*\nu_y^*\otimes End(S_y)\otimes End(V_y)$.  
The linear space $End(S_y)$ can be identified with 
the complex Clifford algebra $\CCl(T_y^*F)$ and then with the complexified exterior 
algebra associated to $T_y^*F$, denoted by $\Lambda^*T_y^*F$. 
Using this identification $k_s(y,y)$ may be considered as a smooth 
function of $s\in\R^+$ with values in 
$\Omega^1_y\otimes \Lambda^*\nu_y^*\otimes \Lambda^*T_y^*F\otimes End(V_y)$. 
Using the riemannian $1$-density  $\mu_y$ at point $y$, this space can be 
identified to $\Lambda^*\nu_y^*\otimes \Lambda^*T_y^*F\otimes End(V_y)$ 
and the corresponding element to $k_s(y,y)$ in this last space is the 
heat kernel, at $(y,y)$ of the rescaled Bismut superconnection as is 
introduced in \cite[page 333]{BeGeVe} without any use of the densities 
(see also remark \ref{ai}). 
Therefore, as long as the $k_s(y,y)$ is concerned, up to the multiplication by 
$\mu_y$ we may assume that our superconnection is the same of \cite[page 333]{BeGeVe}.
In this proof the \emph{vertical degree} of a typical element 
$\omega\otimes\alpha\otimes T$ in 
$\Lambda^*\nu_y^*\otimes \Lambda^*T_y^*F\otimes End(V_y)$
refers to the degree of $\alpha$, its 
\emph{horizontal degree} refers to the degree of $\omega$ 
and the total degree of such element is the sum of its vertical 
and horizontal degrees. 
Let $\delta_s$ be the rescale operator on 
$\Lambda^*\nu_y^*\otimes \Lambda^*T_y^*F\otimes End(V_y)$ that 
multiplies an element of horizontal degree $l$ by $s^{-l}$.
Using the Bismut density theorem for family \cite[theorem 10.21]{BeGeVe}, one has 
the following asymptotic formula when $s$ goes toward $0$ 
\begin{equation}\label{asmoj1}
k_s(y,y)\sim(4\pi s^2)^{-p/2}\sum_{j=0}^\infty s^{2j}\delta_s\,k_j(y) 
\end{equation}
where $k_j(y)$ is a section of 
$\Lambda^*\nu_y^*\otimes \Lambda^*T_y^*F\otimes End(V_y)$ 
such that its total degree is less than or equal to $2j$. 
Moreover the following formula holds, where $F_y$ denotes the curvature of 
the bundle $V$ at point $y$ 
and $\sigma(k_j(y))$ denotes  term of $k_j(y)$ with highest total degree
\begin{equation}\label{asmoj2}
\sum_{j=0}^{p/2}\sigma(k_j(y))=(2\pi i)^{p/2}\operatorname{\hat{A}}(M,F)\,e^{-F_y}~.
\end{equation}
The left hand side of the above formula is called the full symbol of $k_s(y,y)$. 
The action of the supertrace on a typical element 
$\omega\otimes\alpha\otimes T$ equals 
to $str(\alpha)tr(T)\,\omega$. 
It is well known that the supertrace of $\alpha\in \Lambda^qT_y^*F$ as an 
element in $End(S)$ vanishes if $q\neq p$ and equals to its devision by 
the riemannian volume element if $q=p$. Consequently, in above asymptotic 
relation the terms giving rise to negative powers of $s$ have  
no contribution in $str\,k_s(y,y)$ for $j\leq p/2$. 
On the other hand it is clear that the contribution of those terms 
generating a positive power of $s$ vanishes when $t$ goes to $0$. Therefore 
the supertrace $str\,k_s(y,y)$ is coming from constant terms with 
vertical degree $p/2$. These are exactly the terms in the full symbol of $k_s(y,y)$.
So, using the above formula for the full symbole and the above discussion 
one get  
\begin{align*}
\lim_{s\to0}str k_s(y,y)&=(2\pi i)^{p/2}[\operatorname{\hat{A}}(M,F)
\otimes tr e^{-F_y}]_p(y)\\
&=(2\pi i)^{p/2}[\operatorname{\hat{A}}(M,F)\otimes \Ch(E/S)]_p(y)
\end{align*}
where $[~.~]_p$ denotes the part of $.$  
with vertical degree $p$ divided by the volume elements of $T_yF$. Multiplying  
this asymptotic formula by $\mu_y$ is the desired relation \eqref{limmm}.
\end{pf}

\subsection{b-supertrace}\label{bsuper}
With $\tau$ defined by \eqref{tau0}, the expression \eqref{tau} can be
used to define $STr(Q)$ for
$Q$ in the algebra $\mathcal{R}^{-\infty}_{G}
(G^\partial,s^*(\Lambda^*\mu^*)\otimes r^*(\Omega^{\frac{1}{2}}
\otimes E_{|\partial M}))$, i.e. the algebra of the smoothing 
$G$-invariant operators on the boundary groupoid $G^\partial$ 
with rapidly decreasing kernel. 
The proposition \ref{trace1} remains true in this case and 
$\rho\circ STr$ defines a supertrace on this algebra
Due to the presence of the $b$-boundary, 
for $P\in\RR^{-\infty}(G,s^*(\Lambda^*\mu^*)\otimes r^*(\bO\otimes E))$ 
a similar expression in not  
usually convergent and cannot be used to define the trace of $P$. The following 
{\it b-supertrace} is a combination of the Melrose b-trace defect formula 
\cite[proposition 5.9]{Me} 
and the b-supertrace defect formula \cite[proposition 9]{MePi}.
\begin{gather}
\text{}{^b\!}STr:\RR^{-\infty}(G,s^*(\Lambda^*\mu^*)\otimes r^*(\bO\otimes E))
\to C^\infty(M,\Lambda^*\mu^*\otimes\text{}{^b\!}\Omega^1)\notag\\
\text{}{^b\!}STr(P)(x)=\lim_{\epsilon\rightarrow 0}\{\int_{G_x^{t(v)\geq\epsilon}}
\tau(v)\,str\,K_{P|G_x}(v,v)+
\operatorname{ln}\epsilon\int_{G_x^\partial}\tau(v)
\,str\,K_{P|G_x^\partial}(v,v)\}.\label{btrace}
\end{gather}
Here the integrations are taken with respect to longitudinal densities
coming from the supertrace of the Schwartz kernels. So $\text{}{^b\!}STr(P)(x)$
belongs to $\text{}{^b\!}\Omega_x^1(F)$, where $F$ is the
leaf passing by $x$. For keeping the notation simple, we have skipped the 
holonomy action of $v^*$ on the kernel 
(see relation \eqref{tau} for the similar situation with explicit holonomy 
action action ). 
For each $x\in M$, Kernel $\tau\,K_{P|G_x}$ is compactly supported,
so the proof of the lemma 4.62 of \cite{Me} shows that the above
limit exists and  is independent of the boundary defining function
$t$ satisfying
$dt(n)=1$. 
\begin{rem}\label{TrbTr}
If $P$ is an operator as in above with kernel $K_P$ then $\text{}{^b\!}Tr(P)$ 
is defied by relation \eqref{btrace} provided that $str$ is replaced by $tr$. The 
proof of the following proposition remain true if we replace $\text{}{^b\!}STr$ 
by $\text{}{^b\!}Tr$.  
If $K_{P|G_x^\partial}=0$ for each $x\in M$, then the first integral in the 
right hand side of \eqref{btrace} defines $Tr(P)$ provided that 
$str$ is replaced by $tr$. We will use this definition in subsection \ref{dfrabu}.
\end{rem}

\begin{prop}\label{ptrace}
\text{}\\
a) Let $\rho$ be the distributional linear functional on
$C^\infty(M,\text{}{^b\!}\Omega^1\otimes\Lambda^k\nu)$ defined by 
\eqref{gali}. Assume that
the linear functional $\eta$ on $\ob$ defined by
\[\eta(\phi)=\rho(\phi_{|M})\]
is a supertrace on the foliation algebra $\ob$. The linear operator
$\rho\circ\text{}{^b\!}STr$ is independent of $\tau$ and defines a
b-supertrace, i.e. for each $P\in\Diff(G;r^*(\bO\otimes E))$ and for each
$Q\in\RR^{-\infty}(G,r^*(\bO\otimes E))$ the following {\it defect formula} holds
\begin{equation}\label{commutation}
\rho\circ\text{}^bSTr[P, Q]
=\frac{i}{2\pi}\rho\circ\int_\R STr\{\partial_\lambda
I_n(P,\lambda)\circ I_n(Q,\lambda)\}\,d\lambda.
\end{equation}
b) In addition, if $\eta$ is a closed graded trace on $\Omega^k\B$,
then 
\begin{equation}\label{closed}
\rho\circ\text{}^bSTr[\nab,P]=0,
\end{equation}
where $P\in\Omega^{k-1}(r^*\nu^*)\otimes \RR^{-\infty}(G,r^*(\bO\otimes E))$ 
\begin{pf}
The restriction of the kernel $K_P$ of $P\in \RR^{-\infty}$ to the $b$-diagonal of 
$G_b^2$ defines a smooth section $\tilde K_P$ of $End(E)$ over the $b$-diagonal 
of $M_b^2$ such that 
$K_P(v,v)=\tilde K_P(r(v),r(v))$. Using relations
\eqref{char} and \eqref{char1}, the proof of Proposition
\ref{trace1} can be applied to get the following analogue of
\eqref{foot}
\begin{equation}\label{trin}
\rho\circ
\text{}{^b\!}STr(P)=
\lim_{\epsilon\rightarrow 0}\{\int_{M^{t\geq\epsilon}}\rho(x)\,str\tilde K_P(x,x)+
\operatorname{ln}\epsilon\int_{\partial
  M}\rho(y)str\tilde K_P\,(y,y)\}
\end{equation}
The right hand side of this expression  is independent of $\tau$. This 
proves the first part of (a). To prove the defect formula, notice 
that $(PQ)_{|G^\partial}=P_{|G^\partial} Q_{|G^\partial}$ for $P,Q\in
\RR^{-\infty}$. Using the proposition \ref{trace1}, the boundary
part in relation \eqref{btrace} is a trace, i.e. it vanishes on 
supercommutators. In particular this boundary part
has no contribution in the evaluation of $\rho\circ\text{}{^b\!}TR$ on the
supercommutators. Let $\theta\in End(\E)$ be the operator defined by point-wise
multiplication by $\tau$. Using (\ref{btrace}) one has
\[\rho\circ\text{}{^b\!}STr(P)=\rho\circ\text{}{^b\!}STR(\theta P)\]
where $\text{}{^b\!}STR$ is defined by the following relation provided that 
the kernel $K_p$ is compactly supported
\begin{equation}\label{ntrace}
\text{}{^b\!}STR(P)(x)=\lim_{\epsilon\rightarrow
0}{\{\int_{G_x^{t(v)\geq\epsilon}}strK_p(v,v)+
\operatorname{ln}\epsilon\int_{G_x^\partial}strK_p(v,v)\}}
\end{equation}
The second integral in above expression defines the operator
$TR$ on $End_\B(\E_0)$ satisfying 
$\rho\circ Tr\,P_{|G^\partial}=\rho\circ TR(\theta P_{|G^\partial})$. 
As a consequence of the relation (\ref{tau0}) we prove the following equality
\begin{equation}\label{nim}
\text{}{^b\!}STR(\theta QP)=\text{}{^b\!}STR(Q\theta P)
\end{equation}
which implies
\begin{equation}\label{em}
\rho\circ\text{}{^b\!}STr[P,Q]=\rho\circ\text{}{^b\!}STR[\theta P,Q].
\end{equation}
Since the boundary term in the definition of $\rho\circ\text{}{^b\!}STr$ 
has no contribution
in its evaluation on the commutators, the equality \eqref{nim}
is equivalent to the following one
\begin{multline*}
\int_M\rho(x)\int_{G_x^{t\ge\epsilon}}\int_{G_x}du\,dv\tau(u)strK_Q(u,v)\circ
K_P(v,u)\\
=\int_M\rho(y)\int_{G_y^{t\geq\epsilon}}du'\int_{G_y}dv'\tau(v')strK_Q(u',v')\circ
K_P(v',u')
\end{multline*}
 that can be proved as follows
\begin{gather*}
\begin{split}
\int_M\rho(x)&\int_{G_x^{t\ge\epsilon}}\int_{G_x}du\,dv\,\tau(u)trK_Q(u,v)\circ
K_P(v,u)\\
&=\int_M\rho(x)\int_{G_x^{t\geq\epsilon}}\int_{G_x}\int_{G_x}du\,dv\,dw\,\tau(u)\tau(vw^{-1})trK_Q(u,v)\circ
K_P(v,u)\\
&=\int_M\rho(x)\int_{G_x}dw\int_{G_{r(w)}^{t\geq\epsilon}}du'\int_{G_{r(w)}}dv'\,\tau(u'w)\tau(v')
trK_Q(u'w,v'w)\circ K_P(v'w,u'w)\\
&=\int_M\rho(y)\int_{G_y^{t\geq\epsilon}}du'\int_{G_y}dv'\int_{G^y}dw\,\tau(u'w)\tau(v')trK_Q(u',v')\circ
K_P(v',u')\\
&=\int_M\rho(y)\int_{G_y^{t\geq\epsilon}}du'\int_{G_y}dv'\,\tau(v')trK_Q(u',v')\circ
K_P(v',u').
\end{split}
\end{gather*}
It is clear from formula (\ref{ntrace}) that 
\[\text{}{^b\!}STR[\theta P, Q](x)
=\text{}{^b\!}STR_x[(\theta P)_{|G_x}, Q_{|G_x}]~.\]
Because $\tau$ is compactly supported, for each $x\in M$, 
$(\theta P)_{|G_x}$ is a  b-differential operator with compactly supported Schwartz
kernel. It turns out that the proof of
Melrose's {\it trace defect formula} \cite[page 154]{Me} remain true in this compact 
support situation, so the supertrace defect formula \cite[proposition 9]{MePi} 
can be applied to $b$-manifold $G_x$, compactly supported 
smoothing $b$-pseudodifferential operator $\theta P$ and $b$-differential 
operator $Q$ 
to obtain the following leaf-wise defect formula
\[\text{}{^b\!}STR_x[(\theta P)_{|G_x},Q_{|G_x}]
=\frac{i}{2\pi}\int_\R str\,\partial_\lambda
I_n((\theta P)_{|G_x},\lambda).I_n(Q_{|G_x},\lambda)~.\]
Therefore
\begin{align*}
\rho\circ\text{}{^b\!}STr[P,Q]&=\rho\circ \text{}{^b\!}STR[\theta P, Q]\\
&=\frac{i}{2\pi}\int_M\rho(x)\int_\R STr\,\partial_\lambda
I_n(\theta P_{|G_x},\lambda)\circ I_n(Q_{|G_x},\lambda)\,d\lambda\\
&=\frac{i}{2\pi}\int_\R \int_M\rho(x)\,STr\,(\tau(\partial_\lambda
I_n(P,\lambda)\circ I_n(Q,\lambda))_{|G_x^\partial})\,d\lambda\\
&=\frac{i}{2\pi}\int_\R\int_M\rho(x)STr(\partial_\lambda
I_n(P,\lambda)\circ I_n(Q,\lambda))_x\,d\lambda~,
\end{align*}
and this is the desired defect formula \eqref{commutation}.
To prove part (b), it is enough to prove the following analogous of
the relation \eqref{noh}
\begin{equation*}
\rho\circ\text{}{^b\!}STr[\nab,P]=\eta(d^H\circ r^*(\text{}{^b\!}STr(P))).
\end{equation*}
The space $T_v^HG$ is included in $T_vG^\partial$ if $v\in G^\partial$, 
this implies the equality $\nab_{|\partial M}=\nabla^{\E_0}$  
which proves the relation 
\(str\,[\nab,P]_{|G^\partial}=str\,[\nabla^{\E_0},P_{|G^\partial}]\). 
Now, by applying relation \eqref{noh} we get
\begin{equation}\label{vas}
\rho\circ STr[\nab,P]_{G_x^\partial}
=\eta(d^H\circ r^*(STr\,P|G^\partial)).
\end{equation}
This  proves the boundary part of \eqref{closed} which is implicit in the 
definition of $\text{}{^b\!}STr$. However, using
\eqref{char1}, the proof of \eqref{noh} can be slightly modified to get
\begin{equation*}
\rho\circ
STr\,[\nab,P_{|G^{\geq\epsilon}}]=\eta(d^H\circ r^*(STr\,P_{|G^{\geq\epsilon}}).
\end{equation*}
These last two equations and the defining relation \eqref{btrace}
together prove part (b) of the proposition.
\end{pf}
\end{prop}

\section{ Chern character and eta invariant}\label{section4}
\subsection{Chern character}\label{dfrabu}
 Assume that the Dirac operator $D_0$
satisfies the hypothesis \ref{hyp} and fix  $\delta\in\R$ such 
that $0<\delta<\epsilon_0$. With respect to the decomposition $E=E^+\oplus E^-$, 
the parametrix $Q$ given by proposition \ref{frad} is an odd operator of the 
form 
\[Q=\left( \begin{array}{cc}
0&Q^-\\Q^+&0
\end{array}\right), 
\]
In this section, for keeping the notation simple, we use the symbol 
$Q$ to denote $Q^-$.  therefore  
 \begin{gather*}
S_+:=Id-Q\D^+\in\rho_{\ff}\mathcal{R}_{b,G}^{-\infty,
\,\delta}(G,r^*(\bO\otimes E^+))\\
S_-:=Id-\D^+Q\in\rho_{\ff}\mathcal{R}_{b,G}^{-\infty,\,
\delta}(G,r^*(\bO\otimes E^-)).
\end{gather*}
 Since $\delta>0$ we have $I_n(S_\pm,\lambda)=0$ and 
\begin{gather}\label{indo}
I(Q,\lambda)=(i\lambda+D_0)^{-1}.
\end{gather}
Let the projections $p$ and $p_0$ be given by 
\begin{gather*}
p=\begin{pmatrix}S_+^2&S_+(Id+S_+)Q\\S_-\D^+&Id-S_-^2\end{pmatrix}~;~
p_0=\begin{pmatrix}0&0\\0&Id\end{pmatrix}
\end{gather*}
and similar to \cite[relation 5.16]{GoLo-local} define the analytical index class
of the longitudinal Dirac operator by the following formula
\[ind(\Dt):=[p-p_0]
\in K_0(\bar{\mathcal R}_{b,G}).\]
Here $\bar{\mathcal R}_{b,G}$ is the sub algebra of 
${\mathcal R}_{b,G}^{-\infty,\,0}(G,r^*\bO)$ consisting of the operators 
with vanishing indicial family. 
So it is natural  to define the Chern character of the index class by
the following formula where $\mathbb{K}^\E:=(\nab)^2-l$ denotes the curvature of
the connection \nab and $Tr$ is defined in remark \ref{TrbTr}.
\begin{defn}
The Chern character of the index class is defined by
\begin{equation}\label{cern}
\rho\circ\Ch(ind(\D))
=\rho\circ Tr(pe^{-p.\mathbb{K}^\E.p}p-p_0e^{-p_0.\mathbb{K}^\E.p_0}p_0).
\end{equation}
\end{defn}
In this definition the operator which is acted on by $Tr$ is
smoothing with vanishing indicial family, so it has finite $Tr$. 
The index theorem that we want to establish gives a formula to calculate
the Chern character of the index class $ind(\D)$. The 
first step toward this purpose is the following proposition which establishes 
a relationship between the Chern character of the index class and the Chern 
character of the Bismut superconnection given by  
\begin{equation}\label{chersuo}
\rho\circ\bCh(\As)=\rho\circ\str\operatorname e^{-(\As^2-l)}.
\end{equation}
The rest of this subsection is devoted to the proof of this proposition. 
In fact this proposition is a consequence of the next two lemmas. 
\begin{prop}\label{largs}
Let $\rho$ be a functional satisfying the conditions described in 
the proposition (\ref{ptrace}).
The following asymptotic formula holds
\[\lim_{s\rightarrow +\infty}\rho\circ\bCh(\As)
=\rho\circ\Ch(ind(\D)).\]
\end{prop}
To prove this proposition we follow the approach used in 
\cite[section 5]{GoLo-local} and \cite[sections 8-9]{LePi-etale}. Put
\[J=\begin{pmatrix}S_+&-(Id+S_+)Q\\
\D^+&S_-\end{pmatrix}~~;~~p_1=\begin{pmatrix}Id&0\\
0&0\end{pmatrix}\]
we have $p=Jp_1J^{-1}$ where
\[J^{-1}=\begin{pmatrix}S_+&Q(Id+S_-)\\
-\D^+&S_-\end{pmatrix}.\]
Now using the isomorphism $J$ we define the following connection on \E
\[\tnab=(p_1J^{-1}\circ\nab\circ
  Jp_1)+(p_0\nab p_0)
\]
or equivalently
\begin{equation}\label{connec}
J\circ\tnab\circ J^{-1}=p\circ\nab\circ p+(1-p)J\circ\nab\circ
J^{-1}(1-p)
\end{equation}
Put \(\nab=\nab_+\oplus\nab_-\) and \(\tnab=\tnab_+\oplus\tnab_-\), then
\begin{gather}
\tnab_-=p_0\tnab p_0=\nab_-\label{nab1}\\
\tnab_+=p_1\circ\tnab\circ p_1=S_+\nab_+S_++Q(Id+S_-)\nab_-\D^+.\label{nab}
\end{gather}
Finally, using this connection, we introduce the following family of
superconnections from $\E$ into $\Omega^1\E$
\begin{gather*}
A(s)=\begin{pmatrix}\tnab_+&s\D^-\\s\D^+&\tnab_-\end{pmatrix}.
\end{gather*}
Denote its curvature by $K(s)=A^2(s)-l$ and put
\begin{gather*}
e^{-K(s)}=\begin{pmatrix}e^{-K(s)}_{11}&e^{-K(s)}_{12}\\
e^{-K(s)}_{21}&e^{-K(s)}_{22}\end{pmatrix}.
\end{gather*}
With above notation the following formal equalities are easy to
verify
\begin{align*}
STr(e^{-K(s)})&=Tr(e^{-K(s)}_{11})-Tr(e^{-K(s)}_{22})\\
&=Tr(Jp_1e^{-K(s)}p_1J^{-1}-e^{-K(s)}_{22})\\
&=Tr(\D^+e^{-K(s)}_{11}Q(Id+S_-)-e^{-K(s)}_{22}))+Tr(S_+e^{-K(s)}_{11}S_+)
\end{align*}
Although neither $e_{11}^{-K(s)}$ nor $e_{22}^{-K(s)}$ have finite
$Tr$ (they have non-vanishing indicial families),
 the operators in the last line, which are acted on by $Tr$, are rapidly
 decreasing  and have vanishing indicial
operator. Therefore we can
define the Chern character of the superconnection $A(s)$ by the
following formula with value in $C^\infty(M,\Lambda^*\nu^*\otimes\Omega^1)$
\begin{equation}\label{rel}
\Ch(A(s)):=Tr(S_+e^{-K(s)}_{11}S_+)
+  Tr(\D^+e^{-K(s)}_{11}Q(Id+S_-)-e^{-K(s)}_{22})~.
\end{equation}
Our interest in $\Ch(A(s))$ is justified by the following lemma showing  
$\rho\circ \Ch(A(s))=\rho\circ\Ch(ind(\D))$.
\begin{lem}
Under the conditions of the above proposition
\begin{enumerate}
 \item The following relation holds
\[\rho\circ \Ch(A(0))=\rho\circ\Ch(ind(\D)).\]
\item For $s>0$ one has 
\[\rho\circ\Ch(A(s))=\rho\circ \Ch(A(0)).\]
\end{enumerate}
\end{lem}
\begin{pf}
Using the relations (\ref{connec}) and
(\ref{nab}) we have
\[pJe^{-K(0)}J^{-1}p=pe^{-p\mathbb{K}^\E p}p\]
On the other hand $pJ=Jp_1$ and
$p_0e^{-K(0)}p_0=p_0e^{-p_0\mathbb{K}p_0}p_0$ so
\begin{align*}
\Ch(ind(\D))&=
Tr(pe^{-p\mathbb{K}^\E p}p-p_0e^{-p_0\mathbb{K}^\E p_0}p_0)\\
&=  Tr(Jp_1e^{-K(0)}p_1J^{-1}-p_0e^{-K(0)}p_0)\\
&=
Tr(S_+e^{-K(0)}_{11}S_++\D^+e^{-K(s)}_{11}Q(Id+S_-)-e^{-K(s)}_{22})\\
&=\Ch(K(0))
\end{align*}
This proves the first part of the lemma. Now we are going to prove the second part. 
Using the Duhamel formula we get
\begin{align}
\frac{d}{ds} \text{}{^b\!}STr(e^{-K(s)})&=- \text{}{^b\!}STr\int_0^1e^{-uK(s)}
\circ\frac{d}{ds}K(s)\circ e^{-(1-u)K(s)}\,dz\notag\\
&=- \text{}{^b\!}STr\int_0^1[e^{-uK(s)},\frac{d}{ds}K(s)e^{-(1-u)K(s)}]
\,du\label{ho}\\
&~-\text{}{^b\!}STr\frac{d}{ds}K(s).e^{-K(s)}\label{hoho}.
\end{align}
We simplify the expression (\ref{ho}) by means of the defect formula.
>From the very definitions
and up to the smoothing operators with vanishing
indicial families, we have the following formula in the
ungraded notation
\[A^2(s)=\begin{pmatrix}(\tnab_+)^2+s^2\D^-\D^+ &
  s\tnab_+\D^--s\D^-\tnab_-\\0 & s^2\D^+\D^-\end{pmatrix}\]
so
\begin{align*}
I(\frac{d}{ds}K(s),\lambda)&=\frac{d}{ds}I(K(s),\lambda)\\
&=\begin{pmatrix} 2s(\lambda^2+D_0^2) & (\tilde\nabla^{\E_0}-l_+)
(-i\lambda+D_0)-(-i\lambda+D_0)(\tilde\nabla^{\E_0}-l_-)\\
0 &2s(\lambda^2+D_0^2)\end{pmatrix}.
\end{align*}
Using (\ref{nab}) one has also
\begin{equation}\label{yo}
I(e^{-uK(s)},\lambda)=
\begin{pmatrix}
(i\lambda+D_0)^{-1}e^{-u\{(K_0(s)+s^2\lambda+s^2D_0^2)\}}(i\lambda+D_0)
& \mathcal{Z}\\ 0 & e^{-u\{K_0(s)+s^2(\lambda^2+D_0^2)\}}.
\end{pmatrix}
\end{equation}
Using the defect formula and ignoring the
odd functions of $\lambda$, one get
the following equalities
\begin{align}
(\ref{ho}) & =\frac{-i}{2\pi} STr\int_0^1\int_\R\partial_\lambda
           I(e^{-uK_0(s)},\lambda)
           I(\frac{d}{ds}K_0(s)e^{-(1-u)K_0(s)},\lambda)\,d\lambda\,du\notag\\
           &=\frac{-1}{2\pi} Tr\int_\R\int_0^1(i\lambda+D_0)^{-1}
           e^{-u\{K_0+s^2D_0+s^2\lambda^2\}}2s(D_0^2+\lambda^2)
           e^{-(1-u)\{K_0+s^2D_0+s^2\lambda^2\}}\,du\,d\lambda\label{bo}\\
           &+\frac{1}{\pi} Tr\int_\R sD_0
           e^{-\{K_0+s^2D_0+s^2\lambda^2\}}
           \,d\lambda \notag
\end{align}
Apply again the Duhamel and the defect formula to get
\begin{align*}
(\ref{bo})&=-\frac{d}{ds}\text{}{^b\!}Tr[J^{-1},Jp_1e^{-K(s)p_1}]\\
     &=\frac{d}{ds}\text{}{^b\!}STr\,e^{-K(s)}-\frac{d}{ds}\Ch(A(s)),
\end{align*}
where the last equality is coming from relation (\ref{bobo}). Thus
we obtain the equality
\begin{align*}
\frac{d}{ds}\rho\circ \Ch(A(s))=
&-\rho\circ \text{}{^b\!}STr\frac{d}{ds}K(s).e^{-K(s)}\\
&+\frac{1}{\pi}\rho\circ  Tr\int_\R s\tilde D_0
e^{-\{K_0+s^2D_0+s^2\lambda^2\}}\,d\lambda.
\end{align*}
 The following relations are the direct
consequences of the definition of the involved operators and of the
fact that $A(s)$ is a grading reversing operator
\begin{align*}
\text{}{^b\!}STr\frac{d}{ds}K(s).e^{-K(s)}
&=\text{}{^b\!}STr\{(\frac{d}{ds}A(s).A(s)+
A(s)\frac{d}{ds}A(s))e^{-K(s)}\}\\
&= \text{}{^b\!}STr[A(s),\frac{d}{ds}e^{-K(s)}]\\
&= \text{}{^b\!}STr[\tnab-\nab,\begin{pmatrix}0&\D^-\\
\D^+&0\end{pmatrix}]\\
&+ \text{}{^b\!}STr[\nab,\frac{d}{ds}A(s).e^{-K(s)}]\\
&+ \text{}{^b\!}STr[\begin{pmatrix}0&s\D^-\\s\D^+&0\end{pmatrix},
\begin{pmatrix}0&\D^-\\ \D^+&0\end{pmatrix}e^{-K(s)}].
\end{align*}
 The defect formula \eqref{commutation} and the
relation \eqref{yo} imply the vanishing of the first term of the right
hand side of the last equality while the third term is equal to
  \[\frac{1}{\pi}  Tr\int_\R sD_0
    e^{-\{K_0+s^2\tilde D_0+s^2\lambda^2\}}\,d\lambda.\]
Summarizing we get the following relation
\[\frac{d}{ds} \Ch(A(s)= \text{}{^b\!}STr[\nab,\frac{d}{ds}A(s).e^{-K(s)}].\]
 $\rho$ satisfies the condition of the proposition
(\ref{ptrace}) so the relation (\ref{closed}) implies the vanishing
of the above expression when it is acted on by the functional
$\rho$, i.e. 
\[\frac{d}{ds}\rho.\Ch(A(s))=0.\]
This completes the proof of the second part of the lemma
\end{pf}
Since $e^{-K(s)}$ is a smoothing $b$-pseudodifferential operator, 
similar to \eqref{chersuo} 
it is  natural to define the $b$-Chern character of the superconnection 
$A(s)$ by the following formula with value in 
$C^\infty(M,\Lambda^*\nu^*\otimes\text{}{^b\!}\Omega^1)$
\begin{equation*}
\bCh(A(s)):= \text{}{^b\!}STr(e^{-K(s)}).
\end{equation*}
The following lemma shows that the asymptotic behavior of all Chern characters 
defined so far are the same when $s$ goes to $\infty$. Moreover this lemma 
with the previous one provide a proof for the proposition \ref{largs}. 
The following 
relation is easy to verify by using the defect formula for $\text{}{^b\!}Tr$
(see remark \ref{TrbTr}) and will be used in the following discussion
\begin{equation}\label{bobo}
\bCh(A(s))-\Ch(A(s))=  \text{}{^b\!}Tr [J^{-1}, Jp_1e^{-K(s)}p_1].
\end{equation}
\begin{lem}Under the conditions of the above proposition
\begin{enumerate}
\item One has
 \[\rho\circ \bCh(A(s))-\rho\circ \Ch(A(s))=B_1(s)~,\]
where $B_1(s)$ is a boundary depending term going to zero when
$s\rightarrow\infty$.
\item The following relation holds 
\[\rho\circ \bCh(A(s))=\rho\circ \bCh(\As)+B(s)\]
where $B_2(s)$ is a boundary depending term going to zero 
when \marginpar{2}$s\to\infty$.
\end{enumerate}
\end{lem}
\begin{pf}
Using the relations \eqref{cli1} 
and \eqref{cli2},
the operator $K(s)$ has the following expansion in the  ungraded notation
\begin{equation}\label{K(s)}
K(s)=\begin{pmatrix}(\tnab_+)^2-l_++s^2\D^-\D^+&\mathcal{Z}\\
-s\D^+\tnab_++s\tnab_-\D^+&s^2\D^+\D^-+(\tnab_-)^2-l_-\end{pmatrix}
\end{equation}
Using the relations (\ref{nab}) we get
\begin{align*}
K_{21}(s)&=-s\D^+(S_+\nab_+S_++Q(Id+S_-)\nab_-\D^+)+s\tnab_-\D^+\\
         &=-s\D^+S_+\nab_+S_+-s(Id-S_-^2)\nab_-\D^++s\nab_-\D^+\\
         &=-s\D^+S_+\nab_+S_++sS_-^2\nab_-\D^+
\end{align*}
In particular $I(K(s)_{21},\lambda)=0$, so
\begin{align*}
I(e^{-K(s)},\lambda)_{11}&=e^{-I(K(s),\lambda)_{11}}\\
      &=e^{-(\tilde{\nabla}^{\E_0})^2-l_+-s^2(-i\lambda+D_0)(i\lambda+D_0)}\\
                         &=e^{-K_0^-(s)-s^2(\lambda^2+D_0^2)},
\end{align*}
where $K_0^-(s):=(\tilde{\nabla}^{\E_0})^2-l_+$. Now using (\ref{bobo}) and
the defect formula (\ref{commutation}) we obtain
\begin{align*}
\bCh(A(s))-\Ch(A(s))&=  \text{}{^b\!}Tr [J^{-1}, Jp_1e^{-K(s)}p_1]\\
&=\frac{i}{2\pi}
Tr\int_\R\partial_\lambda I(J^{-1},\lambda)\circ I(Jp_1
e^{-K(s)}p_1,\lambda)\,d\lambda\\
&=\frac{1}{2\pi}  Tr\int_\R
e^{-K_0^-(0)-s^2\lambda^2-s^2D_0^2}(i\lambda+D_0)^{-1}\,d\lambda \\
&=\frac{1}{2\pi}Tr\{e^{-K_0^-(0)-s^2D_0^2}\int_\R e^{-s^2\lambda^2}(i\lambda
+D_0)^{-1}\,d\lambda\}
\end{align*}
Denote the last expression by $B(s)$ which is  a boundary depending
term. The self adjoint Dirac operator $D_0$ verifies
the invertibility condition of hypothesis \ref{hyp}. So the $\lambda$-depending
family of operators $(i\lambda+\tilde D_0)^{-1}$ is uniformly
bounded. This implies the finiteness of the integral factor of
$B(s)$.
  In the other hand, using the Volterra development, we have
\begin{gather*}
\begin{split}
e^{-K_0^-(0)-s^2D_0^2}&=
e^{-s^2D_0^2}
-\int_0^1 e^{-u_0 s^2 D_0^2}
(K_0^-(0))e^{-(1-u_0)s^2D_0^2}
+\ldots\\
&+(-1)^k\int_{\triangle_k}e^{-u_0s^2D_0^2}(K_0^-(0))
e^{-u_1s^2D_0^2}(K_0^-(0))\ldots e^{-u_ks^2D_0^2}
\,d\sigma(u_0,\ldots,u_k)+\ldots
\end{split}
\end{gather*}
where $\triangle_k$ denotes the k-simplex
\(\triangle_k=\{(u_0,u_1,\ldots u_k)\in
[0,1]^k~;~~\sum_{j=0}^ku_j=1\}\).\\
We recall from \cite[relation 29]{Lo-h} that the following estimate 
holds for each finite order $G$-invariant differential
operator $P$ 
\begin{equation*}
|P_ue^{-s^2\tilde D^2_0}\tilde D_0^{2l})(u,v)|\leq
 const.F_N(R(u,v))\sum Q(s,s^{-1},\delta)e^{-s^2\delta^2}.
\end{equation*}
Here $Q$ is a polynomial of $s$ and $s^{-1}$ depending on the
differential operator $P$, while $F_N(r)$ is a function of $r$ which
is $O(r^{N+1})$ as $r\rightarrow +\infty$.
Using this estimate, the limit of each term in above expression is
zero when $s\rightarrow\infty$. 
Since $K_0^-(0)$ is a differential two-form, one has
only a finite number of non zero terms in this expression, so  the
vanishing of
all terms implies the vanishing of $B(s)$ when $s$ goes
toward $+\infty$. This completes the proof of the first part of the lemma. 
The second part can be proved by
applying essentially the same computations used  in the first part and in the
previous lemma.
\end{pf}

\subsection{The eta invariant}
Fix $\epsilon>0$ sufficiently small, put $N=[\frac{p}{4}]+1$
and $R(u,v)=d(u,v)-\epsilon$ where $(u,v)\in G^\partial\times_s
G^\partial$. One has the following uniform estimate for this family
which is a consequence of the finite propagation speed estimate explained in
\cite[relation 30]{Lo-h}. This estimate holds only when the hypotheses 
\ref{hyp2} and \ref{hyp} are satisfied. 
\begin{align}
|\nabla_u^\alpha&e^{-s^2D^2_0}D^{2l}_0)(u,v)|\leq\notag\\
&const.(R/s^2)^{-1/2}[R^{-2(k+l)}+R^{-2(k+l)-4N}+\label{asl}\\
&R^{2(k+l)}s^{-4(k+l)}+R^{2(k+l)+4N}s^{-4(k+l)-8N}]e^{-R^2/4s^2}.\notag
\end{align}
If in addition the boundary Dirac operator $D_0$ satisfies the
hypothesis (\ref{hyp}), and if $0<\delta<\epsilon_0$, then the following
estimate holds for each finite order $G$-invariant differential
operator $P$ \cite[relation 29]{Lo-h}
\begin{equation}\label{asl2}
|P_ue^{-s^2D^2_0}D_0^{2l})(u,v)|\leq
 const.F_N(R(u,v))\sum Q(s,s^{-1},\delta)e^{-s^2\delta^2}.
\end{equation}
Here $Q$ is a polynomial of $s$ and $s^{-1}$ depending on the
differential operator $P$, while $F_N(r)$ is a function of $r$ which
is $O(r^{N+1})$ as $r\rightarrow +\infty$.

We need also some estimates for the heat kernel when the parameter
$s$ is small. Using \cite[relation 31]{Lo-h} there is $s_0>0$ and
$\epsilon >0$ such that for $0<s<s_0$ and $d(u,v)<\epsilon$
\begin{equation}\label{asl3}
|\nabla^\alpha_ue^{-s^2D_0}(u,v)|\leq const.s^{-(p+|\alpha|)}
e^{-\frac{d^2(u,v)}{5s^2}}.
\end{equation}
Given the heat kernel of the longitudinal Dirac operator $D_0$, 
the heat kernel for the
superconnection  is given by the Volterra development 
\begin{equation}\label{ik0}
e^{-s^2D_0^2-sF_0-J_0}=e^{-s^2D_0^2}+\sum_{k>0}(-s^2)^kI_k
\end{equation}
where
\begin{multline}\label{ik}
s^{2k}I_k(u,u')=\int_{\triangle_k}\int_{G^k}e^{-s^2t_0
D_0^2}(u,u_0)s^2P_0e^{-s^2t_1D_0^2}(u_0,u_1)s^2P_0\\
\cdots e^{-s^2t_{k-1}D_0^2}s^2P_0e^{-s^2t_kD_0^2}(u_{k-1},u')\,dt.
\end{multline}

 Here $s^2P_0=sF_0+J_0$ and $\triangle_k$ denote the $k$-simplex
\[\triangle_k=\{(t_0,t_1,\ldots t_k)\in
[0,1]^k~;~~\sum_{j=0}^kt_j=1\}.\] If $k>n-p$ then $I_k=0$, so we
have only to prove that for each $s>0$ and for each $0\leq k\leq
n-p$ the  integral (\ref{ik}) is convergent and defines a smooth
section on $(0,+\infty)\times G^\partial\times_s G^\partial$. For a
fixed $k+1$-tuple $(t_0,t_1,\cdots,t_k)$, the rapid decay property
expressed in relation \eqref{asl} implies the absolute convergence
of the integration on $G^k$. In the other hand, for $s\neq0$ at
least one of the heat operator appearing in the integrand is
smoothing (since $t_i\neq0$ for some $0\leq i\leq k$), so the
integrand is smooth with respect to all its variables. Thus for
$s\neq0$,  $I_k$ defines a smoothing operator  which depends
smoothly on $s$.\\
It is clear from above Volterra development that the heat kernel of the rescaled 
$\CCl(1)$-superconnection $\mathbb B_s$ is of the form given by 
relation \eqref{gali}. 
So we can apply the $\CCl(1)$-supertrace \eqref{clsuptr}  
to define the (rescaled) $\Cl(1)$-Chern character by the following relation 
(see \cite[pages 28,33 ]{MePi})
\[\Ch_\alpha(\mathbb B_s):=STr_\alpha(\operatorname e^{-(\mathbb B_s^2-l)}).\] 
The operator 
$\mathbb B_s^2-l$ is even (i.e. grading preserving) on 
$C_c^\infty(G^\partial,\Lambda^ks^*\nu^*\otimes
r^*(\Omega^{\frac{1}{2}} \otimes(E_0\oplus E_0)))$, therefore 
$\Ch_\alpha(\mathbb B_s)$ is a sum of
differential forms of odd orders as we have explained just after the 
relation \eqref{clsuptr}. The following proposition
shows the relation between the Chern character of the
Bismut superconnection and of the \Cl(1)-superconnection (see 
\cite[Proposition 11]{MePi})
\begin{prop}\label{variation}
Let $\rho$ satisfy the condition described in the first part of the
proposition (\ref{ptrace}). We have
\begin{equation*}
\frac{d}{ds}\,\rho\circ \bCh(\As)=-\rho\circ
\text{}{^b\!}STr[\nab, \frac{d\As}{ds}e^{-(\As^2-l)}]
-\frac{1}{2}\rho\circ\hat{\eta}(s)
\end{equation*}
where the rescaled  eta form is defined by
\begin{equation*}
\hat{\eta}(s)=\frac{1}{\sqrt{\pi}} STr_\alpha(\frac{d\mathbb B _s}{ds}
e^{-(\mathbb B _s^2-l_0)}).
\end{equation*}
\begin{pf}
Using Duhamel's formula we get 
\begin{gather}\label{aich}
\frac{d}{ds}\,\rho\circ \bCh(\As)
=-\int_0^1 \rho\circ\text{}{^b\!}STr\{e^{-t(\As^2-l)}
\frac{d(\As^2-l)}{ds} e^{-(1-t)(\As^2-l)}\}\,dt.
\end{gather}
>From \eqref{ind12} the indicial family of the curvature operator
$(\As^2-l)$ equals $\lambda^2+Id.D_0$ which is an even function of
$\lambda$, so the following expression is an odd function of
$\lambda$
\[\partial_\lambda I_n(e^{-t(\As^2-l)},\lambda)\circ
I_n(\frac{d(\As^2-l)}{ds}
e^{-(1-t)(\As^2-l)},\lambda),\]
and the defect formula of the proposition \ref{ptrace} gives the
following expression for (\ref{aich})
\[\frac{d}{ds}\,\rho\circ
\bCh(\As)=-\rho\circ\text{}{^b\!}STr\,\frac{d(\As^2-l)}{ds}
e^{-(\As^2-l)}.\]
$\As$ is an odd operator which commutes with $e^{-(\As^2-l)}$, so we get 
the following relations where
$\A_{[1]}$ is nothing else than \nab:
\begin{align}
\frac{d(\As^2-l)}{ds}e^{-(\As^2-l)}&=\frac{d\As}{ds}e^{-(\As^2-l)}\As+\As
\frac{\As}{ds}e^{-(\As^2-l)}\notag\\
        &=[\As,\frac{d\As}{ds}e^{-(\As^2-l)}]_{gr}\notag\\
        &= [\As-\A_{[1]},\frac{d\As}{ds}e^{-(\As^2-l)}]_{gr}+
[\nab,\frac{d\As}{ds}e^{-(\As^2-l)}]_{gr}\label{moush}
\end{align}
Notice that the operators
\[\As-\A_{[1]}~;~\frac{d\As}{ds}~;~(\As^2-l)\]
are all $G$-invariant b-differential operators. Using again the
relation \eqref{ind12} and  the defect formula, the first term in
the expression \eqref{moush}, followed by operator
$\rho\circ\text{}{^b\!}STr$, is equal to
\begin{align}
\frac{-i}{2\pi}\rho\circ STr&\,\int_\R\partial_\lambda
I(\As-\A_{[1]},\lambda).I(\frac{d\As}{ds}e^{-(\As^2-l)},
\lambda)\,d\lambda\notag\\
&=\frac{-i}{2\pi}\rho\circ STr\int_\R
s\gamma.\frac{d\mathbb B_s}{ds}
e^{-(\mathbb B_s^2-l_0+s^2\lambda^2)}\,d\lambda\notag\\
&=\frac{-i}{2\sqrt{\pi}}\rho\circ
  STr\circ\gamma.\frac{d\mathbb B_s}{ds}e^{-(\mathbb B_s^2-l_0)}\label{davaa}.
\end{align}
The second equality  is a consequence of the fact that the even part
of $\frac{d}{ds}I(\As,\lambda)$, with respect to $\lambda$, is
$\frac{d\mathbb B_s}{ds}$, while the other terms are all odd functions of
$\lambda$. On the other hand
\[str\circ\gamma=tr\circ R\circ\gamma=-i\,str_\alpha,\]
where $R=\begin{pmatrix}1&0\\0 &-1\end{pmatrix}$.  
So $STr\circ\gamma=-i\,Str_\alpha$ and the expression (\ref{davaa}) 
is equal to 
\[\frac{-1}{2\sqrt{\pi}}\rho\circ
  STr_\alpha\frac{d\mathbb B_s}{ds}e^{-(\mathbb B_s^2-l_0)}\]
which is, by definition, equal to
  $\frac{-1}{2}\rho(\hat{\eta}(s))$.
\end{pf}
\end{prop}
\begin{cor}\label{eta}
Let $\rho$ satisfies the conditions described in  part (b) of the proposition
(\ref{ptrace}) then
\[\frac{d}{ds}\,\rho\circ \bCh(\As)
=-\frac{1}{2}\rho(\hat{\eta}(s)).\]
\end{cor}
To establish the index theorem we need to integrate the above
$s$-depending  eta form over $[0,+\infty)$, so the
integrability has to be proved. To prove the next proposition we
follow the proof of \cite[Proposition 25]{Lo-h}.
\begin{prop}
Assume that the $G$-invariant boundary Dirac operator $D_0$ satisfies the
hypothesis \ref{hyp} and that the holonomy groupoid $G$ satisfies 
the polynomial growth conditions of hypothesis \ref{hyp2}.
Then the following integral, taking its value in
$C^\infty(M,\Lambda^*\nu^*)$, is convergent.
\begin{equation}\label{port}
\eta_0=\int_0^{+\infty}\tilde{\eta}(s)\,ds
\end{equation}
\end{prop}
\begin{pf}
Because $\tau$ is compactly supported, it is enough to prove that for 
each $u\in G^\partial$ the following s-depending
 differential form is absolutely integrable on $(0,+\infty)$
\begin{equation}\label{etaa}
str_\alpha(\frac{d\mathbb B_s}{ds}e^{-(\mathbb B_s^2-l_0)})(u,u)
\in \Omega^1(TG^\partial_u)\otimes\Lambda_u^*r^*\nu^*.
\end{equation}
The operator $e^{-(\Bbs^2-l_0)}$ is defined through the Volterra
formula (\ref{ik0}). Using this formula, we consider
a general term in the development of
$\frac{d}{ds}(\mathbb{B}_s)e^{-(\Bbs^2-l_0)})$ and prove that it is
integrable over $[0,+\infty)$. A general term has the following form
 \begin{align*}
s^{2k}\int_{\triangle_k}&(s\alpha.
D_0+\frac{1}{4is^2}\cl(T^H_0))e^{-s^2t_0D_0^2}s^2\\
& P_+e^{-s^2t_1D_0^2}s^2\cdots
P_+e^{-s^2t_{k-1}D_0^2}s^2P_+e^{-s^2t_kD_0^2}(u,u)\,dt_0\cdots dt_k.
\end{align*}
For a fixed $s>0$ take $s^2t_j$ as a new variable and denote it again
by the  symbol $t_j$.
We rewrite the above expression in the following form where $\chi$ is the
characteristic function of the subset $\{0\}\subset\R$
\begin{multline}\label{dob}
\int_{\R^{k+1}}\int_{G_x^\partial}\cdots\int_{G_x^\partial}\chi(s^2-\sum_jt_j)
(s\alpha.D_0+\frac{1}{4is^2}\cl(T^H_0))e^{-t_0D_0^2}(u,u_1)\\
P^+e^{t_1D_0^2}(u_1,u_2)\cdots
P^+ e^{-t_{k-1}D_0^2}(u_{k-1},u_k)P^+
e^{-t_kD_0^2}(u_k,u)\,du_1\cdots\,du_k\,dt_0\cdots dt_k.
\end{multline}
In the relation (\ref{asl3}) take $s_0\leq\frac{s^2}{k+1}$.
Divide the domain of the integration into $2^{k+1}$ pieces
such that on each one the  value of $t_j$ be less than or equal to, or
greater than $s_0$. We begin by the piece on which $t_j>s_0$ for all
 $j$. Using the estimate
(\ref{asl2}) the above integral will be bounded from above by a
constant multiple of
\begin{multline*}
\int_{G_x^\partial}\cdots\int_{G_x^\partial}
F_N(d(u,u_1))F_N(d(u_1,u_2))\cdots F_N(d(u_k,u))\,du_1\,du_2\cdots du_k\\
\times\int_{s_0}\cdots\int_{s_0}\chi(s^2-\sum_jt_j)Q(t_0,t_0^{-1},
\cdots,t_k,t_k^{-1},\delta)
\,e^{-t_0\delta^2}\cdots e^{-t_k\delta^2}\,dt_0\cdots dt_k.
\end{multline*}
$F_N(r)$ is $O(r^{-N})$ when $r\rightarrow +\infty$. So the hypothesis
\ref{hyp2} implies the finiteness of the first factor of the above expression.
The second factor is trivially finite. Moreover it is dominated by
$e^{-s^2\delta^2}$, so for each $u\in G^\partial$, the integration
over this piece provides a function of $s$ which is absolutely integrable 
at $+\infty$.

Now let one of the $t_j$'s, say $t_0$, is less than or equal to $s_0$
and  all others are greater than $s_0$. In this case, using 
estimate (\ref{asl3}), the expression
(\ref{dob}) is bounded from above by
\begin{align}
&\int_0^{s_0}\alpha.D_0\,e^{-t_0D_0^2}(u,u_1)\,dt_0\label{ood}\\
&\times\int_{G_x^\partial}\cdots\int_{G_x^\partial}
F_N(d(u_1,u_2))\cdots F_N(d(u_k,u))\,du_2\cdots du_k\notag\\
&\times\int_{s_0}\cdots\int_{s_0}\chi(s^2-\sum_{j=0}t_j)Q(s,t_1,t_1^{-1},
\cdots,t_k,t_k^{-1},\delta)
\,e^{-t_1\delta^2}\cdots e^{-t_k\delta^2}\,dt_0\cdots dt_k.\notag
\end{align}
To obtain this expression we have used the fact that $\cl(T^H)$
commutes with $D_0^2$ for putting together the operator
$1/s^2\,\cl(T^H)$ and $e^{-t_j\delta^2}$ where $t_j\geq s^2_0/k$.
As before, the second factor is a uniformly bounded function of
$(u_1,u)$, while the third factor is bounded by a constant multiple
of  $e^{-s^2\delta^2/2}$. For $d(u,u_1)\leq2\epsilon$, using
(\ref{asl3}), the first factor can be bounded from above by (see
\cite[relation 42]{Lo-h})
\[const.\int_0^{s_0}t_0^{-(p+1)/2}e^{-\frac{d(u,u_1)^2}{5t_0}}\,dt_0\leq
const.\frac{1}{(d(u,u_1))^{p-1}}\]
In particular, using the polar coordinates of $\R^p$, the left hand 
side of the above expression is integrable
with respect to $u_1$.
If $d(u,u_1)\geq2\epsilon$, using the estimate (\ref{asl}),
the expression (\ref{ood}) is bounded by
\[Q'(R,R^{-1},t_0,t_0^{-1})e^{-R^2/4t_0}\]
where $Q'$ is a polynomial function of its variables. In particular
this expression defines an integrable function of $u_1$ and goes to
zero faster than any power in $d(u,u_1)$ as $d(u,u_1)\rightarrow 0$.
Therefore in this case, for each $u\in G^\partial$, the expression
(\ref{dob}) defines a differential form depending on $s$ which is
integrable at $+\infty$. Clearly this argument can be applied when 
$t_j\geq s_0$ for some $0\leq j\leq k$. But such a $j$ always exists 
because $\sum_jt_j=s^2$ and $s^2\geq(k+1)s_0$. This proves the
integrability at $+\infty$. The following lemma proves the convergence of 
the differential form (\ref{etaa}) when $s\rightarrow0$. This prove the 
integrability at $0$ and completes the proof of the proposition.
\end{pf}
\begin{lem}
The eta form
$\hat{\eta}(s)$ has a limit when $s\rightarrow 0$.
\end{lem}
\begin{pf}
In what follows we use the notation and the content of the proof of  
the proposition \ref{lim}. At first we prove an asymptotic relation
for the $\CCl(1)$-supertrace of the $\CCl(1)$-superconnection $\mathbb B$. 
For this purpos 
consider the foliation $(\partial M\times\R, \partial F\times\R)$ 
which is equipped with the longitudinal product metric $g_0+d^2t$, 
where $t$ denotes 
the coordinate of the $\R$-factor. The connection $T^h$, lifted along 
the $R$-factor, defines a connection for this foliation. 
The vector bundle $E_{|\partial M}\simeq E_0\oplus E_0$ also can be lifted 
along $\R$-factor, and we denote the lifted bundle by the same symbole 
$E_0\oplus E_0$. For $j=1,2$ 
let $(\gamma_j,\beta_j)$ be two paths in a leaf of this foliation 
with common initial and end points. These paths are assumed to be equivalent 
if the holonomy of $\gamma_1$ and $\gamma_2$, 
with respect to the foliation $(M,F)$, are equal. The set of all 
equivalence classes is a groupoid $\tilde G$ over the base manifold $M\times\R$.
We have for example $\bar G_{(x,t)}=G_x^\partial\times\R$ and so on. 
Let $\bar{\mathbb A}$ denote the Bismut superconnection associated to the 
groupoid $\bar G$ and let $\bar {\mathbb A}_s$ denote the associated 
rescaled operator. 
>From the above construction it is clear that the groupoid 
$\bar G$ and the operator $\bar {\mathbb A}_s$ are models for a collar 
neighborhood of $G_c^\partial\in G$ and the restriction of $\As$ to this 
neighborhood. As differential operators on 
$C^\infty(\bar G, r^*(\Omega^{1/2}\otimes\Lambda^*\nu^*\otimes E))$ we have  
\[\bar{\mathbb A}_s^2=\mathbb B_s^2-s^2\frac{\partial^2}{\partial t^2}~.\]
Because all involved operators in above commute with $\partial/\partial t$ we have 
\[e^{-(\bar{\mathbb A}_s^2-l)}((u,t),(u,t))
=e^{-(\mathbb B_s^2-l)}(u,u)\,e^{s^2(\partial)^2}(t,t)
=r^*(e^{-\tilde{\mathbb B}_s^2}(y,y))\,e^{s^2(\partial)^2}(t,t)\] 
Since the above expression is local when $t$ goes toward $0$, 
we can and will assume that $(M\times\R,F\times\R)$ is a fibration and that  
$E$ is the lifting of $(S_0\oplus S_0)\otimes V$ along $\R$. 
Therefore using asymptotic relation \eqref{asmoj2} and the fact 
that $e^{s^2(\partial_t)^2}(t,t)=(4\pi s^2)^{-1/2}$ we get 
\begin{equation*}
\tilde k_s(y,y)\sim(4\pi s^2)^{-(p-1)/2}\sum_{j=0}^\infty s^{2j}\delta_s\tilde k_j(y) 
\end{equation*}
Here $\tilde k_s(y,y)$, the heat kernel of $\tilde{\mathbb B}_s^2$, and 
$\tilde k_j(y)$'s are elements of 
$\Omega^1\otimes\Lambda^*\nu^*\otimes End(S_0\oplus S_0)\otimes End(V_y)$. 
Moreover, using the definition \eqref{csup} of $\mathbb B$ and Duhamel's 
formula, it turns out that the $End(S_0\oplus S_0)\otimes End(V_y)$-factor 
has the following form 
\[(\Id R+\alpha\,S)\otimes T=\left(\begin{array}{cc}R&S\\S&R\end{array}\right)
\otimes T\]
where $R~,S\in End(S_0)$ and $T\in End(V_y)$.
By definition the $\CCl(1)$-supertrace (denoted by $str_\alpha$) 
of such an element is 
\[str_\alpha((\Id R+\alpha\,S)\otimes T):=2tr(S).tr(T)~.\] 
With respect to the identification 
$End(S_0\oplus E_0)\simeq\Lambda^*T^*_y\partial F\otimes T^*_t\R$ 
the elements of the form $\Id R$ corespond to 
$(\Lambda^*T^*_y\partial F)\wedge\partial_t$ while the elements of the form 
$\alpha\,S$ correspond to $\Lambda^*T^*_y\partial F$. 
Therefore the $\CCl(1)$-supertrace of an element in 
$\Lambda^*T^*_y\partial F\otimes T^*_t\R$ is non zero only if this 
element is the volume element of $\Lambda^*T^*_y\partial F$ and in this 
case its $\CCl(1)$-supertrace is equal to $1$. Consequently, applying 
the $str_\alpha$ on both side 
of the above asymptotic relation and using the fact that the total degree of 
$k_j(y)$ is a most $2j$, we get the following asymptotic 
formula when $s$ goes to $0$ 
\begin{equation}\label{cloc}
str_\alpha\, e^{-(\mathbb B_s^2-l)}(u,u)
\sim r^*(B_0+sB_1+s^2B_j+\dots)
\end{equation}
Here $B_j$'s are elements in the finite dimensional linear 
space $\Omega^1(T_y\partial F)\otimes\Lambda^*\nu_y^*$ and above asymptotic 
relation occurs in this space.\\
Let $\R^+$ denote the positive real line and consider the foliation
$(\partial M\times \R^+, \partial F)$. Let $t$ denote the
parameter of $\R^+$, then $\bar{g}(x,t):=t^{-1}g_0(x)$ define a
longitudinal Riemannian metric for this foliated manifold
where $g_0$ is the longitudinal Riemannian structure of 
$(\partial M,\partial F)$. Moreover this foliation is equipped with the 
horizontal distribution $T^h\partial M$ and the lifted (along $\R^+$) 
Clifford bundle $E$. The holonomy groupoid associated to 
this foliation is $G\times\R^+$. The $\CCl(1)$-superconnection 
$\mathbb B$ associated to this foliation satisfies the asymptotic  
formula \eqref{cloc}. Using this asymptotic formula instead of the
Bismut local index formula in the proof of
\cite[theorem 10.32]{BeGeVe}, we get the convergence of the eta form when $s$ 
goes toward $0$.
\end{pf}
The above proposition justifies the following definition
\begin{defn}\label{etad}
The   eta invariant, depending on the boundary Dirac
operator $D_0$, the boundary horizontal distribution
$T^h_{|\partial M}$ and the functional $\rho$ is defined by the
following relation
\[\rho(\eta_0)=\rho\circ\int_0^{+\infty}\tilde{\eta}(s)\,ds\]
\end{defn}
\section{The index theorem}\label{section5}
Now we can formulate the main theorem of the paper which is a
direct conclusion of the relation (\ref{evaa}), the
corollary \ref{eta} and the propositions \ref{lim} and \ref{largs}.
\begin{thm}[ b-index theorem for $G$-invariant Dirac operators]\label{hyo}
Let $G$ denote the holonomy groupoid associated to the boundary
foliated manifold $(M,F)$ satisfying the hypothesis \ref{hyp2}.
Assume that $p=\dim F$ is an even number and let $E$ be a
longitudinal Clifford bundle over $(M,F)$ with
 the associated  Dirac operator $\Dt$.
Let $\D=r^*\Dt$ be the $G$-invariant Dirac operator, acting on the
smooth sections of the $G$-invariant vector bundle $r^*E\rightarrow
G$. We assume that the boundary family satisfies the invertibility
hypothesis \ref{hyp}.
 Let $\rho$ be a functional satisfying the conditions described in the
proposition \ref{ptrace}. Then the following index formula holds
\begin{equation*}
\rho(\Ch(ind(\D)))= \frac{1}{(2\pi
i)^{p/2}}\int_M\rho(x)\{\operatorname{\hat{A}}(M,F)
\operatorname{Ch}(E/S)(x)\}_p-\frac{1}{2}\rho(\eta_0).
\end{equation*}
In this formula the Chern character of the analytical index
$\Ch(ind(\D))$ is defined by (\ref{cern}) while the
 eta form $\eta_0$ is given by (\ref{port}). The integrand
is defined by relations \eqref{akar} and \eqref{twi} and the operation 
$\{~.~\}_p$ is defined just before the proposition \ref{lim}.
\end{thm}
\begin{rems}\text{}\\
i)When $\partial M=\emptyset$ and under the polynomial growth condition 
of hypothesis \ref{hyp2}, this index theorem implies the
index theorem established by A.Gorokhovsky and J.Lott in
\cite[Theorem 1]{GoLo}. The extra coefficient $\frac{1}{(2\pi
i)^{p/2}}$ appeared in our formula is a consequence of our
definition for the Chern character (Compare the relation \eqref{a0}
in the appendix with the
relation (34) of \cite{GoLo}).\\
ii) Let $\mathcal{C}^0$ be an invariant zero current, i.e. an holonomy
invariant measure. Using this invariant measure, the relation
\eqref{majid} defines a functional $\rho: C_0(M, \Omega^1)\rightarrow
\R$ satisfying the conditions of the proposition \ref{ptrace}. In this
case, again with hypothesis \ref{hyp2}, the above theorem reduces 
to the A.Connes' foliation index
theorem \cite{Cos} when $\partial M=\emptyset$.
\end{rems}
\section{Appendix}\label{fram}
In this appendix we describe briefly a general construction for the
Chern character which is proposed by Gorokhovsky and Lott
\cite{GoLo}. This is a generalization of the classical Chern-Weil
construction.

Let $(\mathcal{B}^*,d)$ be a $\Z$-graded algebra where $d$ is a
linear operator of grading one and $d^2\neq 0$. We denote the
operator $d^2$ by $\alpha$. The algebra $\mathcal{B}$ can be viewed
as a
 $\Z_2$-graded algebra  with the even sub algebra
$\mathcal{B}^e=\oplus \mathcal{B}^{2i}$ 
and the odd sub space 
$\mathcal{B}^o=\oplus \mathcal{B}^{2i+1}$. 
It is clear that $d$ is a grading reversing 
operator on $\mathcal{B}$.
 Let $\mathcal{E}$ be a $\Z_2$-graded left $\mathcal{B}$-module equipped with a
superconnection $\nabla$. This means that
\(\nabla:\mathcal{E}\rightarrow\mathcal{E}\) is an odd operator satisfying the
following relation for all $\phi\in\mathcal{B}$ and $\xi\in\mathcal{E}$
\[\nabla(\phi\xi)=d(\phi)\xi+\phi\nabla(\xi).\]
The algebra ${\rm End}_\mathcal{B}(\mathcal{E})$ is naturally a
$\Z_2$-graded algebra and is equipped with a canonical supertrace
\begin{gather*}
str:{\rm End}_{\mathcal{B}}(\mathcal{E})\rightarrow
\mathcal{B}_{ab}\\
str(\xi^*\otimes \xi)=\xi^*(\xi)~;~\xi^*\in\E^*,~\xi\in\E
\end{gather*}
where $\mathcal{B}_{ab}:=\B/[\B,\B]$.

In the classical case, when $\alpha=0$, one defines the Chern
character of this superconnection by the following relation
\begin{equation}\label{a0}
\Ch_\rho(\mathcal{E},\nabla)=\rho\circ str(e^{-K})\in \mathcal A\tag{A.1}
\end{equation}
where $K:=\nabla^2\in {\rm End}_{\mathcal{B}}(\mathcal{E})$ denotes
the curvature of the superconnection $\nabla$ and  $\rho$ is an
appropriate graded trace on $\mathcal{B}$ which takes its values in
an Abelian group $\mathcal{A}$. A graded trace means a linear
function which vanishes on the supercommutators as well as on the
image of $\alpha$. When $\alpha\neq 0$ we assume the existence of an
even linear map $l:\mathcal{E}\rightarrow \mathcal{E}$ satisfying
the following relations
\begin{gather}
l(\phi\xi)=\alpha(\phi)\xi+\phi l(\xi)\tag{A.2}\\
l\circ\nabla=\nabla\circ l\tag{A.3}
\end{gather}
and  define the curvature of the connection
$\nabla$ to be $K=\nabla^2-l$. The following lemma shows that $K$
shares the common properties with the classical curvature 
(for the proof see lemmas 2 and 7 of \cite{GoLo}).
\begin{lem}\label{kop}
\text{}\\
a) $K:\mathcal{E}\rightarrow\mathcal{E}$ is $\mathcal{B}$-linear.\\
b) $K$ is a flat $\B$-linear map, i.e $[\nabla,~K]=0$.
\end{lem}
We can now  define the Chern character by the formula (\ref{a0}).
Let $d^t$ denote the adjoint of $d$.
 Given a graded trace $\rho$, it is clear
that $d^t(\rho):=\rho\circ d$ is a  graded trace too, moreover
$d^t\circ d^t=0$. Thus the set
of all graded traces form a complex with differential $d^t$. The
trace $\rho$ is closed if $d^t(\rho)=0$. For the proof of the following 
lemma we refere to \cite[lemma 8]{GoLo}.
\begin{lem}
\text{}\\
a)For a fixed closed graded trace $\rho$, the Chern
character defined in (\ref{a0}) is independent of the connection
$\nabla$, so one can denote the Chern character by $\Ch_\rho(\mathcal{E})$.\\
b)If $\rho_1$ is homologous with $\rho_2$ 
then \(\Ch_{\rho_1}(\mathcal{E})=\Ch_{\rho_2}(\mathcal{E})\).
\end{lem}

\end{document}